\documentclass[12pt]{amsart}
\usepackage{mathrsfs}
\usepackage{amsfonts}
\usepackage{amssymb}
\usepackage[papersize={8.2in,11in},textwidth=16.4cm,textheight=22cm,centering]{geometry}
\usepackage{enumerate}

\usepackage[colorlinks=true,citecolor=blue,linkcolor=blue]{hyperref}
\hypersetup{
pdfstartpage=1,
pdfstartview=FitH}

\usepackage{color}

\DeclareFontFamily{U}{mathx}{\hyphenchar\font45}
\DeclareFontShape{U}{mathx}{m}{n}{
      <5> <6> <7> <8> <9> <10>
      <10.95> <12> <14.4> <17.28> <20.74> <24.88>
      mathx10
      }{}
\DeclareSymbolFont{mathx}{U}{mathx}{m}{n}
\DeclareMathAccent{\widecheck}{\mathalpha}{mathx}{"71}

 \usepackage{caption} 
\numberwithin{equation}{section}

\allowdisplaybreaks

\newtheorem*{theorem*}{Theorem}
\newtheorem{theorem}{Theorem}[section]
\newtheorem{lemma}{Lemma}[section]

\newtheorem{definition}{Definition}[section]
\newtheorem{proposition}{Proposition}[section]

\newtheorem{conjecture}{Conjecture}[section]

\newtheorem*{claim*}{Claim}

\theoremstyle{remark}
\newtheorem{remark}{\bf Remark}
\newcommand{\ud}{\mathrm{d}}
\newcommand{\ue}{\mathrm{e}}

\newcommand{\tr}{\mathrm{tr}}
\newcommand{\sym}{\mathrm{sym}}

\newcommand{\kl}{\mathrm{Kl}}
\newcommand{\Frob}{\mathrm{Frob}}

\DeclareMathOperator*{\Res}{Res}
\DeclareMathOperator{\Mod}{mod}

\renewcommand{\bmod}[1]{\,(\Mod{ #1})}

\newcommand{\bk}{\mathbf{k}}

\newcommand{\by}{\mathbf{y}}

\newcommand{\bC}{\mathbf{C}}
\newcommand{\bF}{\mathbf{F}}

\newcommand{\bQ}{\mathbf{Q}}
\newcommand{\bR}{\mathbf{R}}
\newcommand{\bZ}{\mathbf{Z}}

\newcommand{\cA}{\mathcal{A}}
\newcommand{\cB}{\mathcal{B}}
\newcommand{\cC}{\mathcal{C}}
\newcommand{\cD}{\mathcal{D}}
\newcommand{\cE}{\mathcal{E}}
\newcommand{\cF}{\mathcal{F}}
\newcommand{\cG}{\mathcal{G}}
\newcommand{\cH}{\mathcal{H}}

\newcommand{\cL}{\mathcal{L}}
\newcommand{\cM}{\mathcal{M}}
\newcommand{\cN}{\mathcal{N}}
\newcommand{\cP}{\mathcal{P}}

\newcommand{\cR}{\mathcal{R}}
\newcommand{\cS}{\mathcal{S}}

\newcommand{\fc}{\mathfrak{c}}
\newcommand{\ff}{\mathfrak{f}}

\newcommand{\fq}{\mathfrak{q}}

\newcommand{\fA}{\mathfrak{A}}

\newcommand{\fS}{\mathfrak{S}}

\begin{document}

\title{When Kloosterman sums meet Hecke eigenvalues}
\author{Ping Xi}

\address{Department of Mathematics, Xi'an Jiaotong University, Xi'an 710049, P. R. China}\email{ping.xi@xjtu.edu.cn}


\subjclass[2010]{11L05, 11F30, 11N36, 11T23}

\keywords{Kloosterman sum, Hecke--Maass eigenvalue, equidistribution, Selberg sieve, Riemann Hypothesis over finite fields}

\begin{abstract}
By elaborating a two-dimensional Selberg sieve with asymptotics and equidistributions of Kloosterman sums from $\ell$-adic cohomology, as well as a Bombieri--Vinogradov type mean value theorem for Kloosterman sums in arithmetic progressions, it is proved that for any given primitive Hecke--Maass cusp form of trivial nebentypus, the eigenvalue of the $n$-th Hecke operator does not coincide with the Kloosterman sum $\kl(1,n)$ for infinitely many squarefree $n$ with at most $100$ prime factors. This provides a partial negative answer to a problem of Katz on modular structures of Kloosterman sums.

\end{abstract}

\dedicatory{{\it \small Dedicated to Professor \'Etienne Fouvry \\ on the occasion of his sixty-fifth birthday }}

\maketitle


\section{Introduction}\label{sec:introduction}

We are concerned with the {\it normalized} Kloosterman sum
\[\kl(a,c)=\frac{1}{\sqrt{c}}~\sideset{}{^*}\sum_{u\bmod c}\ue\Big(\frac{au+\overline{u}}{c}\Big)\]
defined for all $c\in\bZ^+$ and $a\in\bZ.$ Denote by $\cP$ the set of primes. For each $p\in\cP$ and $a\in\bZ$, the celebrated Weil's bound asserts that $|\kl(a,p)|\leqslant2,$ from which one finds there exists a certain $\theta_p(a)\in[0,\pi]$ such that 
\[\kl(a,p)=2\cos\theta_p(a).\]
In his famous lecture notes, Katz \cite[Chapter 1]{Ka80} proposed the following three problems (with $a\neq0$ fixed):
\begin{enumerate}[(I)]
\item Does the density of $\{p\in\cP:\kl(a,p)>0\}$ in $\cP$ exist? If yes, is it equal to $1/2?$

\item Is there a measure on $[0,\pi]$ such that $\{\theta_p(a):p\in\cP\}$ equidistributes?

\item Consider the Euler product
\[L_a(s):=\prod_{p\in\cP,p\nmid a}\Big(1-\frac{\kl(a,p)}{p^{s}}+\frac{1}{p^{2s}}\Big)^{-1}\]
for $\Re s>1.$ Is it defined to be an $L$-function attached to some Maass form of level $\fq$ with $\fq$ being a power of 2?
\end{enumerate}

Problem I is also referred to the sign change problem of Kloosterman sums, and the first serious progress was made by Fouvry and Michel \cite{FM03b,FM07}, who proved that there are at least $\gg X/\log X$ squarefree numbers $c\in[X,2X]$ with $\omega(c)\leqslant23$ such that $\kl(1,c)>0$ (resp. $\kl(1,c)<0$), where $\omega(c)$ denotes the number of distinct prime factors of $c$. The method of Fouvry and Michel includes a pioneer combination of the Selberg sieve, spectral theory of automorphic forms and $\ell$-adic cohomology.
The constant 23 was later sharpened by Sivak-Fischler \cite{SF09}, Matom\"aki \cite{Ma11} and the author \cite{Xi15}, and the current record 7 is due to the author \cite{Xi18}. Quite recently, Drappeau and Maynard \cite{DM19} reduced the constant further to 2 by assuming the existence of Landau--Siegel zeros in a suitable way.

Problem II concerns the horizontal equidistribution of Kloosterman sums, and it is expected the Sato--Tate measure $\ud\mu_{\text{ST}}:=\frac{2}{\pi}\cos^2\theta\ud\theta$ does this job. In fact, Katz \cite[Conjecture 1.2.5]{Ka80} formulated a precise conjecture that for each fixed integer $a\neq0$, the set $\{\theta_p(a):p\in\cP,p\nmid a\}$ of Kloosterman sum angles should equidistribute with respect to $\ud\mu_{\text{ST}}$. It then follows immediately from this conjecture that the desired density in Problem I is $1/2$; i.e.,
\[\lim_{x\rightarrow+\infty}\frac{|\{p\in\cP\cap[1,x]:\kl(a,p)>0\}|}{|\cP\cap[1,x]|}=\frac{1}{2}.\]

The original Sato--Tate conjecture was first formulated independently by Sato
and Tate in the context of elliptic curves, and then reformulated and extended to
the framework of Hecke eigencuspforms for $SL_2(\bZ)$ by Serre \cite{Se68}, predicting the similar equidistributions 
of Fourier coefficients of such cusp forms.
Very recently, the conjecture has been confirmed by L. Clozel, M. Harris \& R. Taylor \cite{CHT08} for non-CM elliptic curves over $\bQ$ with non-integral $j$-invariants, and was later generalized by
Barnet-Lamb, Geraghty, Harris and Taylor \cite{BGHT11} for non-CM, holomorphic elliptic modular newforms of weight $k\geqslant2$, level $N$. Much earlier before this resolution, the vertical analogue with $p$ fixed and the form varying was considered independently
by Conrey, Duke and Farmer \cite{CDF97} and Serre \cite{Se97}. In parallel with the vertical Sato--Tate distribution for cusp forms, Katz \cite{Ka88} proved for Kloosterman sums that the set $\{\theta_p(a):1\leqslant a<p\}$ becomes equidistributed with respect to $\ud\mu_{\text{ST}}$ as long as $p\rightarrow+\infty.$

In view of the similarity between the distributions of Kloosterman sums and Hecke eigenvalues of holomorphic cusp forms, it seems natural to expect that $\pm \kl(1,p)$ might coincide with the $p$-th Fourier coefficient of some holomorphic Hecke cusp form. In fact, thanks to Deligne \cite{De80}, $-\kl(1,p)$ and eigenvalue $\lambda_f(p)$ of the $p$-th Hecke operator acting on a primitive holomorphic cusp form $f$ are both known to be Frobenius traces of an $\ell$-adic Galois representations of rank 2 and weight 0. Unfortunately, it is easily known that they could not coincide, since $\kl(1,p)$ cannot lie in any fixed number field (see \cite{Bo00} for some discussions).
Problem III of Katz concerns the modular structure of Kloosterman sums, and predicts that the situation might be valid if one considers Maass forms in place of holomorphic ones. 

In what follows, we take $f$ to be a primitive Hecke--Maass cusp form of level $\fq,$ trivial nebentypus and eigenvalue $\lambda=1/4+t^2,$ so that it is a joint eigenfunction of the Laplacian and Hecke operators. Suppose $f$ admits the following Fourier expansion
\begin{align*}f(z)=\sqrt{y}\sum_{n\neq0}\lambda_f(n)K_{it}(2\pi|n|y)\ue(nx),\end{align*}
where $\lambda_f(1)=1$ and $K_\nu$ is the $K$-Bessel function of order $\nu$. The trivial nebentypus enables $\lambda_f$'s to be real numbers.
As eigenvalues of Hecke operators, $\lambda_f$'s are expected to satisfy the inequality
\begin{align}\label{eq:RamanujanPetersson}
|\lambda_f(n)|\leqslant n^\vartheta\tau(n)
\end{align}
for some $\vartheta<1/2.$ The Ramanujan--Petersson conjecture asserts that $\vartheta=0$ is admissible, and the current record, due to Kim--Sarnak \cite{KS03}, 
takes $\vartheta=7/64.$ Although it is already known that most Hecke--Maass cusp forms $f$ satisfiy \eqref{eq:RamanujanPetersson} with $\vartheta=0$ (see Sarnak \cite{Sa87}), the distribution of $\lambda_f(n)$ is still mysterious in many aspects. Problem III is thus two-fold: $\lambda_f(n)$ is suggested to be controlled
by virtue of Kloosterman sums; and conversely, spectral theory of Maass forms might be helpful to understand {\it non-trivial} analytic information about the Euler product $L_a(s),$
which would yield non-trivial progresses towards to Problems I and II.
Unfortunately, Problem III seems too optimistic to be true, but there seems no satisfactory approach that has been found to attack it; and to my best knowledge, the only known result regarding this problem was obtained by Booker \cite{Bo00} based on numerical computations:  if $\kl(1,p)=\pm\lambda_f(p)$ for some primitive Hecke--Maass cusp form $f$ of level $\fq=2^\nu$ and eigenvalue $\lambda,$
then $\fq\cdot(\lambda+3)>18.3\times10^6.$

In this paper, we present an analytic-number-theoretic approach to Problem III, which enables us to provide a partial negative answer with almost primes in place of primes.

\begin{theorem}\label{thm:non-identity}
Let $f$ be a primitive Hecke--Maass cusp form $f$ with trivial nebentypus.
Then there exist infinitely many squarefree number $n$ with at most $100$ prime factors, such that 
\begin{align*}\lambda_f(n)\neq \pm\kl(1,n).\end{align*}
Quantitatively, for $\eta\in\{-1,1\},$ there exists certain constant $c=c(f)>0,$ such that
\begin{align*}|\{n\in[X,2X]:\lambda_f(n)>\eta\cdot\kl(1,n),~\omega(n)\leqslant 100,~\mu^2(n)=1\}|\geqslant \frac{cX}{\log X}\end{align*}
and
\begin{align*}|\{n\in[X,2X]:\lambda_f(n)<\eta\cdot\kl(1,n),~\omega(n)\leqslant 100,~\mu^2(n)=1\}|\geqslant \frac{cX}{\log X}\end{align*}
hold for all $X>1/c.$
\end{theorem}

In fact, we can prove the following general theorem.
\begin{theorem}\label{thm:non-identity-general}
For any $\eta\in\bR$ and each primitive Hecke--Maass cusp form $f$ of trivial nebentypus,
there exist two constants $c=c(f,\eta)>0$ and $r=r(\eta)<+\infty,$ such that
\begin{align*}|\{n\in[X,2X]:\lambda_f(n)>\eta\cdot\kl(1,n),~\omega(n)\leqslant r,~\mu^2(n)=1\}|\geqslant \frac{cX}{\log X}\end{align*}
and
\begin{align*}|\{n\in[X,2X]:\lambda_f(n)<\eta\cdot\kl(1,n),~\omega(n)\leqslant r,~\mu^2(n)=1\}|\geqslant \frac{cX}{\log X}\end{align*}
hold for all $X>1/c.$ In particular, one may take $r(\pm\frac{1}{2018})=25,$ $r(\pm2018)=41.$
\end{theorem}

Theorems \ref{thm:non-identity} and \ref{thm:non-identity-general} are new, at least to the author, even if there is no restriction on the size of $\omega(n).$
The merit of Theorem \ref{thm:non-identity-general} is revealed by the flexibility of $\eta$. Although we cannot provide a complete negative answer to Problem III, it seems that there is little hope to find a suitable Hecke--Maass cusp form to capture modular structures of Kloosterman sums following the line in Problem III. However, the function field analogue was confirmed by 
Chai and Li \cite{CL03} that the relevant Kloosterman sums (defined over the residue field of completion of the function field at place $v$) and Hecke eigenvalues of a certain $GL_2$ automorphic form can coincide up to a negative sign.

In a private communication, Katz proposed a problem to consider an analogue of Problem III with the cubic exponential sum
\[B(a,c):=\frac{1}{\sqrt{c}}\sum_{x\bmod c}\ue\Big(\frac{x^3+ax}{c}\Big)\] 
in place of $\kl(1,p).$ The vertical Sato--Tate distribution of $B(a,p)$, as $a$ runs over $(\bZ/p\bZ)^\times$ for sufficiently large prime $p$, was already proved by Katz \cite[Section 7.11]{Ka90}, and it is also expected that the horizontal equidistribution of $B(a,p)$ should be true. However, to prove analogues of Theorems \ref{thm:non-identity} and \ref{thm:non-identity-general} seems beyond our current approach. 

In fact, the proofs of Theorems \ref{thm:non-identity} and  \ref{thm:non-identity-general} rely on a kind of Bombieri--Vinogradov type equidistribution 
for Kloosterman sums $\kl(1,c)$ (see Lemma \ref{lm:BV-FM} below), and this was proved by Fouvry and Michel \cite{FM07} by appealing to the spectral theory of automorphic forms. It is thus natural to expect such a theorem should also exist for $B(1,c)$. We would like to mention a similar result due to Louvel \cite{Lo14} that such Bombieri--Vinogradov type equidistribution
holds for cubic exponential sums modulo Eisenstein integers, for which he employed the spectral theory of cubic metaplectic forms, and cubic residue symbols can be well-introduced. However, as Louvel has pointed out, it is not yet known how to move from the cubic exponential sums modulo Eisenstein integers to those modulo rational integers in the horizontal aspect.

Theorems \ref{thm:non-identity} and \ref{thm:non-identity-general}  will be proved by 
appealing to a weighted Selberg sieve and the arguments will be outlined in the next section.

\subsection*{Notation} As usual, $\ue(z)=\ue^{2\pi iz}$ and $\mu,\varphi,\tau$ denote the M\"obius, Euler and divisor functions, respectively. We use $\omega(n)$ to count the number of distinct prime factors of $n$. The superscript $*$ in summation indicates primitive elements.
Given $X\geqslant2$, we set $\mathcal{L}=\log X$ and the notation $n\sim N$ means $N<n\leqslant2N.$  
For a sequence of coefficients $\boldsymbol\alpha=(\alpha_m)$, denote by $\|\cdot\|_1$ and $\|\cdot\|$ the $\ell_1$- and $\ell_2$-norms, respectively, i.e., $\|\boldsymbol\alpha\|_1=\sum_m|\alpha_m|,\|\boldsymbol\alpha\|=(\sum_m|\alpha_m|^2)^{1/2}.$

\subsection*{Acknowledgements} 
I am very grateful to Professors \'Etienne Fouvry, Nicholas Katz and Philippe Michel for their valuable suggestions, comments and encouragement. Sincere thanks are also due to the referee for his/her patient comments  and corrections that lead to a much more polished version of this article.
This work is supported in part by NSFC (No. 11601413, No. 11971370).

\smallskip

\section{Setting-up: outline of the proof}\label{sec:outline}
\subsection{A weighted Selberg sieve}
Suppose $(a_n)_{n\leqslant x}$ is a sequence of non-negative numbers. The sieve method was originally designed to capture how often these numbers are supported on primes, although current status only allows us to detect {\it almost primes} in most cases.
A convenient approach was invented by Selberg \cite{Se71} in 1950's in connection with the twin prime conjecture. Precisely, he suggests to consider the weighted average
\begin{align*}
\sum_{n\leqslant x}a_nw_n\{\rho-\tau(n)\},
\end{align*}
where $w_n$ is a non-negative function, and $\rho$ is to be chosen appropriately such that the total average is positive for all sufficiently large $x$, from which one obtains the existence of $n$ such that $\omega(n)\leqslant\log\rho/\log2.$ The ingenuity then lies in the choice of $w_n$, which should attenuate the contributions from those $n$'s that have many prime factors. A typical choice for $w_n$, due to Selberg himself, is the square of the M\"obius transform of a certain smooth truncation of the M\"obius function; see \cite[Chapters 4--7, 10]{HR74} and \cite[Chapter 7]{FI10} for detailed discussions.

Focusing on Problem I of Katz in the first section on sign changes of Klooterman sums, the author \cite{Xi15,Xi18} introduced the above weighted Selberg sieve to the context of Kloosterman sums, in which situation $\tau(n)$ is replaced by a certain truncated divisor function that suits well for the application of Sato--Tate distribution of Kloosterman sums in the vertical aspect. Such experiences motivate us to consider Problem III of Katz in a similar manner.

We need to make some preparations.
\begin{itemize}
\item Let $n$ be a positive integer. For $\alpha,\beta>0$ and $\varDelta>1$, define a truncated divisor function
\begin{align}\label{eq:tau(n;alpha,beta)}
\tau_{\varDelta}(n;\alpha,\beta)=\sum_{\substack{d\mid n\\d\leqslant n^\frac{1}{1+\varDelta}}}\alpha^{\omega(d)}\beta^{\omega(n/d)}.\end{align}

\item Let $X$ be a large number and define $\vartheta\in~]0,\frac{1}{4}]$ by $\sqrt{D}=X^\vartheta\exp(-\sqrt{\cL})$. We choose $(\varrho_d)$ such that
\begin{align}\label{eq:varrho_d}
\varrho_d=
\begin{cases}
\mu(d)\Big(\dfrac{\log(\sqrt{D}/d)}{\log\sqrt{D}}\Big)^2,\ \ & d\leqslant\sqrt{D},\\
0, &d>\sqrt{D}.
\end{cases}
\end{align}

\item Let $\varPsi$ be a fixed non-negative smooth function supported in $[1,2]$ with the normalization
\begin{align}\label{eq:normalization-Psi}
\int_{\bR}\varPsi(x)\ud x=1.
\end{align}
The Mellin transform of $\varPsi$ is defined as
\begin{align*}\widetilde{\varPsi}(s)=\int_{\bR}\varPsi(x)x^{s-1}\ud x.\end{align*}
Hence $\widetilde{\varPsi}(1)=1$. Integrating by parts, we have
\begin{align*}\widetilde{\varPsi}(s)\ll(|s|+1)^{-A}\end{align*}
for any $A\geqslant0$ with an implied constant depending only on $A$ and $\varPsi$. 

\item For any fixed $\eta\in\bR,$ put
\begin{align}\psi(n)=\psi_{f,\eta}(n):=\lambda_f(n)-\eta\cdot \kl(1,n).\end{align}

\item For all $z\geqslant2$, define
\begin{align*}P(z)=\prod_{p<z,p\in\cP}p.\end{align*}
We will specialize $z$ later as a small power of $X$ such that $z^{12}\leqslant X$.
\end{itemize}

Our theorems will be concluded by effective evaluations of the following weighted average
\begin{align}\label{eq:Hpm(X)}
H^\pm(X)=\sum_{n\geqslant1}\varPsi\Big(\frac{n}{X}\Big)\mu^2(n)\{|\psi(n)|\pm\psi(n)\}\{\rho-\tau_{\varDelta}(n;\alpha,\beta)\}\Big(\sum_{d|(n,P(z))}\varrho_d\Big)^2,\end{align}
where $\rho,\vartheta,\alpha,\beta,z,\varDelta$ are some parameters to be chosen later.  Clearly, we have
\begin{align}\label{eq:Hpm(X)-lowerbound}
H^\pm(X)\geqslant\rho\cdot H_1(X)-2H_2(X)\pm\rho\cdot H_3(X)
\end{align}
with
\begin{align*}H_1(X)&=\sum_{n\geqslant1}\varPsi\Big(\frac{n}{X}\Big)\mu^2(n)|\psi(n)|\Big(\sum_{d|(n,P(z))}\varrho_d\Big)^2,\\
H_2(X)&=\sum_{n\geqslant1}\varPsi\Big(\frac{n}{X}\Big)\mu^2(n)|\psi(n)|\tau_{\varDelta}(n;\alpha,\beta)\Big(\sum_{d|(n,P(z))}\varrho_d\Big)^2,\\
H_3(X)&=\sum_{n\geqslant1}\varPsi\Big(\frac{n}{X}\Big)\mu^2(n)\psi(n)\Big(\sum_{d|(n,P(z))}\varrho_d\Big)^2.\end{align*}
Note that $\eta$ is contained implicitly in all above and subsequent formulations; we will keep it fixed and not display this until necessary.
The task reduces to find a positive lower bound for $H_1(X)$ and an upper bound for $H_2(X)$, which should be of the same order of magnitude, and also a reasonable estimate for $H_3(X)$. In fact, we may prove the following three propositions.

\begin{proposition}\label{prop:H1(X)-lowerbound}
For large $X$, we have
\begin{align*}
H_1(X)&\geqslant
(1+o(1))\sum_{2\leqslant i\leqslant7}I_i\cdot\sqrt{l_i^3/{u_i}} \cdot \frac{X}{\log X},
\end{align*}
where $I_i,l_i,u_i$ are given as in Proposition $\ref{prop:Sigma(X,alphai)-lowerbound}.$
\end{proposition}

\begin{proposition}\label{prop:H2(X)-upperbound}
Let $\alpha=\frac{3\pi }{8}$ and $\beta=\frac{1}{2}.$ For large $X$, we have
\begin{align*}H_2(X)&\leqslant (1+o(1))(1+|\eta|/2)\fS\Big(\vartheta,\frac{2\vartheta\log X}{\log z}\Big)\frac{X}{\log X},\end{align*}
where $\fS(\cdot,\cdot)$ is defined by $\eqref{eq:fS(gamma,tau)}$.
\end{proposition}

\begin{proposition}\label{prop:H3(X)-estimate}
For large $X$, we have
\begin{align*}H_3(X)&\ll(1+|\eta|)X(\log X)^{-A}\end{align*}
for any $A>0,$ provided that $\vartheta\leqslant\frac{1}{4},$ where the implied constant depends on $A,f$ and $\varPsi.$ 
\end{proposition}

Upon suitable choices of $\rho,\vartheta$ and $z$, the positivity of $H^\pm(X)$ would imply, for $X$ large enough, that there exists $n\in[X,2X]$ with
\[\tau_{\varDelta}(n;\alpha,\beta)<\rho\]
for which $\lambda_f(n)-\eta\cdot \kl(1,n)>0$ (or $<0$). In Section \ref{sec:numerical}, we will do necessary numerical computations that lead to Theorems \ref{thm:non-identity} and \ref{thm:non-identity-general}.

\subsection{Ingredients of the proof}
The proofs of Propositions \ref{prop:H1(X)-lowerbound}, \ref{prop:H2(X)-upperbound} and \ref{prop:H3(X)-estimate} form the heart of this paper.
The proof of Proposition \ref{prop:H3(X)-estimate} is not new and was stated as \cite[Proposition 2.1]{FM07} in a slightly different setting.
For the proof of Proposition \ref{prop:H1(X)-lowerbound}, we will be restricted to some specialized integers having fixed number of prime factors, for which we may explore the vertical Sato--Tate distribution for Kloosterman sums and moments of Hecke eigenvalues to produce a positive lower bound for $H_1(X)$. 

There is another new ingredient in this paper that the lower bound for $H_1(X)$ relies on the economic control of the correlation
\begin{align*}
\sum_{n\in\cS}\lambda_f(n)\kl(1,n),
\end{align*}
where $\cS$ is a suitable set of the products of a fixed number of distinct primes. It is expected that $\lambda_f(n)$ does not correlate with $\kl(1,n)$ as $n$ runs over $\cS$, and an upper bound which beats the trivial estimate $O(|\cS|)$ is highly desirable. Unfortunately, we do not know how to capture such cancellations, even if $n$ is relaxed to run over consecutive integers. Alternatively, it could be a courageous choice to ignore the sign changes of 
summands, and a suitable upper bound for 
\begin{align*}
\sum_{n\in\cS}|\lambda_f(n)\kl(1,n)|
\end{align*}
with a small scalar
might also suffice.
In fact, for $n$ being a product of distinct primes, say $n=qr$ with $q,r\in\cP$, we may decompose the summand as $|\lambda_f(q)\lambda_f(r)\kl(\overline{r}^2,q)\kl(\overline{q}^2;r)|$. Our observation lies in the fact that 
$|\lambda_f(p)|$ and $|\kl(\overline{p}^2,q)|$ are both smaller than 1, say $\delta$, on average while $p$ runs over a suitable set of primes; see Lemma \ref{lm:pi-kappa(X)-moments} and Lemma \ref{lm:Kloosterman-Mellin-prime} for details.
The factor $\delta^j$ for some large $j$ will then appear if $n$ has more prime factors, and $\delta^j$ can be considerably small if $j$ is taken to be reasonably large. This, while $n$ is restricted to be products of a large number of distinct primes, will lead to quite a small scalar in the upper bound for the above average with absolute values, although we cannot save anything in the order of magnitude. 
Typically, we require $n$ to have 7 distinct prime factors, but this would arise a combinatorial disaster while evaluating $\varrho_d$ in the sieve weight to conclude Proposition \ref{prop:H1(X)-lowerbound}. Thus, the restriction $n\mid P(z)$ in \eqref{eq:Hpm(X)} is introduced to overcome such difficulty in computations.
More precisely, we may restrict $n$ to be the product of certain primes of prescribed sizes larger than $z$, then only $d=1$ survives in the convolution $\sum_{d\mid (n,P(z))}\varrho_d.$

The upper bound for $H_2(X)$ also relies on the vertical Sato--Tate distribution for Kloosterman sums, and by appealing to an idea in our previous work \cite{Xi18}, we may reduce the {\it dimension} of the sifting problem by introducing $\tau_{\varDelta}(n;\alpha,\beta)$ with appropriate choices for $\alpha$ and $\beta$, so that the upper bound for $H_2(X)$ can be controlled more effectively. Due to the appearance of $n\mid P(z)$, one has to evaluate the $k$-dimensional sifting average
\[\sum_{n}\mu^2(n)b_n\Big(\sum_{d\mid (n,P(z))}\varrho_d\Big)^2\]
with $\{b_n\}$ being a non-negative multiplicative function mimicking $k^{\omega(n)}$ on average.
In particular, we may develop an asymptotic evaluation in the case $k=2$ upon the choice \eqref{eq:varrho_d}, which we call a two-dimensional Selberg sieve with asymptotics.
A complete and precise statement will be included as the appendix.

\subsection{Correlations of Kloosterman sums and Hecke eigenvalues}
Before closing this section, we would like to formulate two conjectures which illustrate the correlations between Kloosterman sums and Hecke eigenvalues. 
\begin{conjecture}\label{conj:primes}
Let $f$ be a fixed primitive cusp form $($holomorphic for Maass$).$ For all large $X,$ we have
\begin{align*}
\sum_{p\leqslant X}\lambda_f(p)\kl(1,p)=o(X\cL^{-1}).
\end{align*}
\end{conjecture}
If Conjecture \ref{conj:primes} could be proved affirmatively, it would follow that there exist $100\%$ primes $p$ such that $\lambda_f(p)\neq\kl(1,p)$ for each primitive cusp form $f$, which provides a negative answer to Problem III of Katz.

In order to consider the correlations along prime variables, it should be natural at first to study the average over consecutive integers as we have mentioned. To this end, we may formulate the following correlation with a precise saving.
\begin{conjecture}\label{conj:integers}
Let $f$ be a fixed primitive cusp form $($holomorphic for Maass$).$ For all large $X,$ we have
\begin{align*}
\sum_{n\leqslant X}\lambda_f(n)\kl(1,n)=O(X\cL^{-2018}).
\end{align*}
\end{conjecture}

It seems that the above two conjectures are both beyond the current approach, and the resolutions should require new creations both from automorphic forms and algebraic geometry.

\smallskip
\section{Maass forms and Kloosterman sums}\label{sec:Maass-Kloosterman}

\subsection{Maass forms}
We will not need too much information on Maass forms. The following moments of Fourier coefficients at prime arguments would be most of what is required.

Let $f$ be a primitive Hecke--Maass cusp form $f$ of level $\fq,$ trivial nebentypus and eigenvalue $\lambda$ as an eigenfunction of the Laplacian operator.
For each $\kappa\geqslant0$ and $X>1,$ define
\begin{align*}
\pi_\kappa(X)=\sum_{p\leqslant X}|\lambda_f(p)|^\kappa.
\end{align*}

\begin{lemma}\label{lm:pi-kappa(X)-moments}
For all large $X,$ we have
\begin{align*}
\pi_\kappa(X)=c_\kappa(1+o(1))X\cL^{-1}
\end{align*}
for $\kappa=0,2,4,6$ with $c_0=c_2=1,$ $c_4=2,$ $c_6=5,$ and
\begin{align*}
\pi_\kappa(X)\leqslant c_\kappa(1+o(1))X\cL^{-1}
\end{align*}
for $\kappa=1,3$ with $c_1=\frac{11}{12},$ $c_3=\sqrt{5}.$
\end{lemma}

\proof
We only consider the cases $\kappa\geqslant1.$ Following the approach of Hadamard--de la Vall\'ee-Poussin to the classical prime number theorem, it suffices to consider the non-vanishing and holomorphy of the symmetric power $L$-functions $L(\sym^\kappa f,s)$ at $\Re s=1$ with $\kappa=2,4,6.$
These are already known due to a series of celebrated works \cite{GJ78,KSh00,KSh02,Sh89}. 
In fact, $[0,\pi]$ is identified with the set of conjugacy classes of the compact group $SU_2(\bC)$ via the map
$g\in SU_2(\bC)\mapsto \tr(g)=2\cos\theta;$
the image of the probability Haar measure of $SU_2(\bC)$ is just the Sato--Tate measure $\mu_{\mathrm{ST}}.$
For $\kappa=2j$ $(j=1,2,3)$, we have
$$c_\kappa=\int(2\cos\theta)^{2j}\ud\mu_{\mathrm{ST}}=\frac{1}{j+1}\binom{2j}{j}.$$
In particular, $c_2=1,$ $c_4=2$ and $c_6=5.$

The value of $c_1$ follows from the asymptotics for $\pi_\kappa(X)$ with $\kappa=0,2,4,6$ and the inequality
\begin{align*}
|y|\leqslant\frac{1}{36}(13+29y^2-7y^4+y^6)
\end{align*}
upon the choice by Holowinsky \cite{Ho09},
which is valid for all $y\in\bR$. The value of $c_3$ follows from Cauchy's inequality and the asymptotics for $\pi_6(X).$
\endproof

\subsection{Kloosterman sums}
Following Deligne \cite{De80} and Katz \cite{Ka88}, it is known that
$$a\mapsto-\kl(a,p)=-2\cos\theta_p(a),\ \ a\in \bF_p^\times$$
is the trace function of an $\ell$-adic sheaf $\mathcal{K}l$ on $\mathbf{G}_{m}(\bF_p)=\bF_p^\times$, which is of rank 2 and pure of weight 0. 
Alternatively, we may write
\begin{align*}2\cos\theta_p(a)=\tr(\Frob_a,\mathcal{K}l),\ \ a\in \bF_p^\times.\end{align*}
By Weyl's criterion and the Peter--Weyl
theorem, Katz's vertical equidistribution, as mentioned in the first section,
reduces to control the cancellations within the averages
\begin{align*}\sum_{a\in\bF_p^\times}\sym_k(\theta_p(a))=\sum_{a\in\bF_p^\times}\tr(\Frob_a,\sym^k\mathcal{K}l),\end{align*}
where $\sym^k\mathcal{K}l$ is the $k$-th symmetric power
of the Kloosterman sheaf $\mathcal{K}l$ (i.e., the composition of the sheaf $\mathcal{K}l$ with the $k$-th symmetric power representation of $SL_2$) and
\begin{align*}\sym_k(\theta)=\frac{\sin(k+1)\theta}{\sin\theta}.\end{align*}
In fact, Katz \cite[Example 13.6]{Ka88} proved that
\begin{align}\label{eq:Katz}
\left|\sum_{a\in\bF_p^\times}\sym_k(\theta_p(a))\right|\leqslant\frac{1}{2}(k+1)\sqrt{p}.
\end{align}
It is natural to expect that the square-root cancellation also holds if replacing $\theta_p(a)$ by $\theta_p(\Pi(a))$ for any non-constant rational function $\Pi$ of fixed degree over $\bF_p^\times.$ In general, we have the following estimate.
\begin{lemma}\label{lm:Michel}
Let $\psi$ and $\chi$ be additive and multiplicative characters $($not necessarily non-trivial$)$ modulo $p$ and $\Pi$ a fixed non-constant rational function modulo $p$. For each fixed positive integer $k,$ there exists some constant $B$ depending only on $\deg(\Pi),$ such that
\begin{align}\label{eq:psi-Kl}
\sideset{}{^*}\sum_{a\bmod p}\psi(a)\sym_k(\theta_p(\Pi(a)))\ll k^B\sqrt{p}
\end{align}
\begin{align}\label{eq:chi-Kl}
\sideset{}{^*}\sum_{a\bmod p}\chi(a)\sym_k(\theta_p(\Pi(a)))\ll k^B\sqrt{p}
\end{align}
hold with implied constants depending at most on $B.$ In particular, one can take $B=1$ if $\Pi(a)=1/a^2.$
\end{lemma}

The case $\Pi(a)=1/a^2$ in \eqref{eq:psi-Kl} was contained in Michel \cite[Corollarie 2.9]{Mi98b} and there is no essential difference when extending to general $\Pi.$ The bound in Lemma \ref{lm:Michel} lies in the fact that the underlying sheaf
$\sym^k([\Pi^*\mathcal{K}l])$ is of rank $k+1$, while the Artin--Scherier sheaf $\cL_\psi$ is of rank 1 if $\psi$ is non-trivial. 
These two geometrically irreducible sheaves are not geometrically isomorphic, and the square-root cancellation then follows from the 
Riemann Hypothesis of Deligne \cite{De80} (see also \cite[Theorem 4.1]{FKM15}, for instance, for practical use in analytic number theory).
The bound \eqref{eq:chi-Kl} follows by noting that the Kummer sheaf $\cL_\chi$ is geometrically irreducible and of rank 1 if $\chi$ is non-trivial.

For $(a,c)=1$, define
\begin{align}\label{eq:Omega}
\varOmega(a,c):=\kl(\overline{a}^2,c).\end{align}
It follows from the Chinese remainder theorem that the twisted multiplicativity $\varOmega(a,rs)=\varOmega(ar,s)\varOmega(as,r)$ holds for all $a,r,s$ with $(r,s)=(a,rs)=1.$ For each Dirichlet character $\chi\bmod c$, define
\begin{align*}
\widetilde{\varOmega}(\chi,c):=\frac{1}{\sqrt{c}}~\sideset{}{^*}\sum_{r\bmod c}\overline{\chi}(r)|\varOmega(r,c)|.
\end{align*}
For $\chi_1\bmod{c_1}$ and $\chi_2\bmod{c_2}$ with $(c_1,c_2)=1,$
the Chinese remainder theorem yields
\begin{align}\label{eq:Kl-Mellin-multiplicativity}
\widetilde{\varOmega}(\chi_1\chi_2,c_1c_2)=\chi_1(c_2)\chi_2(c_1)\widetilde{\varOmega}(\chi_1,c_1)\widetilde{\varOmega}(\chi_2,c_2).
\end{align}
For prime moduli, we have the following asymptotic characterizations.

\begin{lemma}\label{lm:Kloosterman-Mellin-prime}
Let $p$ be a large prime. Then
\begin{align*}
\widetilde{\varOmega}(\chi,p)=\delta_\chi \sqrt{p}+O(\log p),
\end{align*}
where $\delta_\chi$ vanishes unless $\chi$ is the trivial character mod $p,$ in which case it is equal to $\frac{8}{3\pi},$ and the implied constant is absolute.
\end{lemma}

\proof In view of Lemma \ref{lm:Michel}, we may apply Lemma \ref{lm:cos-Chebyshev} with
\[J=\varphi(p), \ \ B=1,\ \ U=c\sqrt{p},\]
\[\{y_j\}_{1\leqslant j\leqslant J}=\{\overline{\chi}(r):1\leqslant r\leqslant p-1\}, \ \ \{\theta_j\}_{1\leqslant j\leqslant J}=\{\theta_p(\overline{r}^2):1\leqslant r\leqslant p-1\},\]
where $c$ is a large suitable constant,
so that
\begin{align*}
\sqrt{p}\widetilde{\varOmega}(\chi,p)-\frac{8}{3\pi}~\sideset{}{^*}\sum_{r\bmod p}\overline{\chi}(r)
&\ll\sqrt{p}\log K+\frac{p^{3/2}}{K}
\end{align*}
for any $K>1$.
The proof is completed by taking $K=p.$
\endproof

\begin{lemma}\label{lm:Kloosterman-Mellin-composite}
Let $q\geqslant2$ be a squarefree number and $\chi$ a primitive character mod $q$. Then we have
\begin{align*}
|\widetilde{\varOmega}(\chi,q)|\leqslant c^{\omega(q)}\log q
\end{align*}
for some absolute constant $c>0.$
\end{lemma}

\proof
In view of \eqref{eq:Kl-Mellin-multiplicativity}, we have
\begin{align*}
|\widetilde{\varOmega}(\chi,q)|=\prod_{p\mid q}|\widetilde{\varOmega}(\chi_p,p)|,
\end{align*}
where $\chi_p$ is a non-trivial character mod $p$. By Lemma \ref{lm:Kloosterman-Mellin-prime}, there exists some absolute constant $c_0>0$, such that
\begin{align*}
|\widetilde{\varOmega}(\chi,q)|
\leqslant\prod_{p\mid q}c_0\log p=c_0^{\omega(q)}\prod_{p\mid q}\log p
\leqslant c_0^{\omega(q)}\sum_{d\mid q}\mu^2(d)\log d
\leqslant(2c_0)^{\omega(q)}\log q.
\end{align*}
This completes the proof of the lemma by taking $c=2c_0$.
\endproof

The following bilinear form estimation can be found in \cite[Corollaire 2.11]{Mi95} and a more general statement has been proved in \cite[Theorem 1.17]{FKM14}.

\begin{lemma}\label{lm:bilinear}
Let $p$ be a large prime and $(a,p)=1.$
For each $k\geqslant1$ and any coefficients $\boldsymbol\alpha=(\alpha_m),\boldsymbol\beta=(\beta_n),$ we have
\begin{align*}\mathop{\sum_{m\sim M}\sum_{n\sim N}}_{(mn,p)=1}
\alpha_m\beta_n\sym_k(\theta_p(\overline{(amn)^2}))\ll 
\|\boldsymbol\alpha\|\|\boldsymbol\beta\|(MN)^{\frac{1}{2}}(p^{-\frac{1}{4}}+N^{-\frac{1}{2}}+M^{-\frac{1}{2}}p^{\frac{1}{4}}(\log p)^{\frac{1}{2}}),\end{align*}
where the implied constant depends polynomially on $k$.
\end{lemma}

\begin{remark}
Lemma \ref{lm:bilinear} is non-trivial as long as $N>\log p,M>p^{\frac{1}{2}}(\log p)^2$ and $p>\log(MN).$
\end{remark}

The following lemma is originally proved by Fouvry and Michel \cite[Proposition 7.2]{FM07} using $\ell$-adic cohomology.
\begin{lemma}\label{lm:bilinearform-twisted}
Suppose $q=q_1q_2\cdots q_s$ with $q_1,q_2,\cdots,q_s$ being distinct primes. For each $s$-tuple of positive integers $\bk=(k_1,k_2,\cdots,k_s),$ and any coefficients $\boldsymbol\alpha=(\alpha_m),\boldsymbol\beta=(\beta_n),\boldsymbol\gamma=(\gamma_{m,n})$ with $m\equiv m'\bmod n\Rightarrow\gamma_{m,n}=\gamma_{m',n},$ we have
\begin{align*}\mathop{\sum_{m\sim M}\sum_{n\sim N}}_{(mn,q)=1}\alpha_m\beta_n\gamma_{m,n}&\prod_{1\leqslant j\leqslant s}\sym_{k_j}(\theta_{q_j}(\overline{(mnq/q_j)^2}))\\
&\ll c(s;\mathbf{k})
\|\boldsymbol\alpha\|\|\boldsymbol\beta\|\|\boldsymbol\gamma\|_{\infty}(MN)^{\frac{1}{2}}(q^{-\frac{1}{8}}+N^{-\frac{1}{4}}q^{\frac{1}{8}}+M^{-\frac{1}{2}}N^{\frac{1}{2}}),\end{align*}
where $c(s;\mathbf{k})=3^s\prod_{j=1}^s(k_j+1)$ and the implied constant is absolute.
\end{lemma}

\begin{remark}
Lemma \ref{lm:bilinearform-twisted} is non-trivial as long as $M>N\log q>q^{\frac{1}{2}}(\log q)^2$ and $q>\log(MN).$
\end{remark}

\begin{lemma}\label{lm:bilinearform-avergeovermoduli}
Let $P,M\geqslant3.$ Suppose $\boldsymbol\gamma=(\gamma_p)$ is a complex coefficient
 supported on primes in $]P,2P]$ and $\Pi$ is a fixed non-constant rational function with integral coefficients in numerators and denominators. Then there exists some constant $B=B(\deg(\Pi))>0,$ such that for each $k\geqslant1$ and arbitrary coefficient $\boldsymbol\alpha=(\alpha_m)$ supported in $]M,2M],$ 
\begin{align*}
\sum_{p\sim P}\gamma_p\sum_{m\sim M}\alpha_m \sym_k(\cos\theta_p(\Pi(m)))\ll k^B(M^{\frac{1}{2}}+P\log P)\|\boldsymbol\alpha\|\|\boldsymbol\gamma\|
\end{align*}
holds with some implied constant depending at most on $B.$
\end{lemma}

\begin{remark}
A typical situation is $\gamma_p\equiv1$, in which case Lemma \ref{lm:bilinearform-avergeovermoduli} becomes non-trivial as long as $P,M/(P\log^2P)\rightarrow+\infty.$ It is an important and challenging problem to beat the barrier $M=P$ for a general coefficient $\boldsymbol\alpha=(\alpha_m)$. We would like to mention a deep result of Michel \cite{Mi98a}, who considered the special case $k=1,$ $\boldsymbol\gamma\equiv 1,$ $\Pi(m)=m$, and he was able to work non-trivially even when $M$ is quite close to $\sqrt{P}.$
\end{remark}

\proof
Write $K(m,p)=\sym_k(\cos\theta_p(\Pi(m)))$ and
denote by $S$ the average in question. First, by Cauchy's inequality, we have
\begin{align}\label{eq:Cauchy}
|S|^2\leqslant \|\boldsymbol\alpha\|^2\varSigma
,\end{align}
where
\begin{align*}
\varSigma=\sum_{m\sim M}\left|\sum_{p\sim P}\gamma_pK(m,p)\right|^2.\end{align*}
Squaring out and switching summations, we find
\begin{align*}
\varSigma=\mathop{\sum\sum}_{p_1,p_2\sim P}\gamma_{p_1}\overline{\gamma}_{p_2}\sum_{m\sim M}K(m,p_1)\overline{K(m,p_2)}=\varSigma^=+\varSigma^{\neq},\end{align*}
where we split the double sum over $p_1,p_2$ according to $p_1=p_2$ or $p_1\neq p_2.$

Trivially, we have
\begin{align}\label{eq:Sigma=}
\varSigma^==\sum_{p\sim P}|\gamma_p|^2\sum_{m\sim M}|K(m,p)|^2\leqslant (k+1)^2M\|\boldsymbol\gamma\|^2.
\end{align}
By completion method, we may derive, for $p_1\neq p_2,$ that
\begin{align*}
\sum_{m\sim M}K(m,p_1)\overline{K(m,p_2)}
&=\sum_{r\bmod{p_1p_2}}K(r,p_1)\overline{K(r,p_2)}\sum_{\substack{m\sim M\\ m\equiv r\bmod{p_1p_2}}}1\\
&=\frac{1}{\sqrt{p_1p_2}}\sum_{|h|\leqslant \frac{1}{2}p_1p_2}\sum_{m\sim M}\ue\Big(\frac{hm}{p_1p_2}\Big)\widehat{K}(h\overline{p_2},p_1)\overline{\widehat{K}(-h\overline{p_1},p_2)},\end{align*}
where
\begin{align*}
\widehat{K}(y,p)=\frac{1}{\sqrt{p}}~\sideset{}{^*}\sum_{x\bmod p}K(x,p)\ue\Big(\frac{-xy}{p}\Big).\end{align*}
From Lemma \ref{lm:Michel} it follows that
\begin{align*}
\varSigma^{\neq}
&\leqslant k^B\mathop{\sum\sum}_{p_1\neq p_2\sim P}
\frac{|\gamma_{p_1}\gamma_{p_2}|}{\sqrt{p_1p_2}}\sum_{|h|\leqslant \frac{1}{2}p_1p_2}\min\Big\{M,\frac{p_1p_2}{h}\Big\}\\
&\ll k^B(M+P^2\log P)\|\boldsymbol\gamma\|^2.\end{align*}
Combining this with \eqref{eq:Sigma=}, we find
\begin{align*}
\varSigma\ll k^B(M+P^2\log P)\|\boldsymbol\gamma\|^2,\end{align*}
from which and \eqref{eq:Cauchy}, the lemma follows immediately.
\endproof

\begin{lemma}\label{lm:bilinearaverageoverprimes}
Let $P,X\geqslant3.$ Suppose $\boldsymbol\gamma=(\gamma_p)$ is a complex coefficient
supported on primes in $]P,2P]$ and $\nu$ is a multiplicative function such that
\[\sum_{n\leqslant N}\tau(n)|\nu(n)|^2\ll N(\log N)^\kappa\]
for some constant $\kappa\geqslant1.$  Then we have
\begin{align*}
\sum_{p\sim P}\gamma_p\sum_{\substack{n\sim X\\(n,p)=1}}\mu^2(n)&\nu(n)\Lambda(n)|\varOmega(n,p)|=
\frac{8}{3\pi}\sum_{p\sim P}\gamma_p\sum_{\substack{n\sim X\\(n,p)=1}}\mu^2(n)\nu(n)\Lambda(n)\\
&\ \ \ \ +O\Big(\cL^{A}\{PX^{\frac{1}{2}}+P^{\frac{1}{4}}X+P^{\frac{1}{2}}X\cL^{-2A}+(PX)^{\frac{3}{4}}\}\|\boldsymbol\gamma\|\Big)
\end{align*}
for any $A>\kappa+2$, where the implied constant depends only on $A$ and $\kappa$.
\end{lemma}

\begin{remark}
Lemma \ref{lm:bilinearaverageoverprimes} is non-trivial as long as $\cL\ll P\ll X\cL^{-3A}.$
\end{remark}

\proof
In view of the Chebyshev approximation for $|\cos|$ (see Lemma \ref{lm:cos-Chebyshev}), it suffices to consider
\begin{align*}
\sum_{p\sim P}\gamma_p\sum_{\substack{n\sim X\\(n,p)=1}}\mu^2(n)&\nu(n)\Lambda(n)\sym_k(\cos\theta_p(\overline{n^2})).
\end{align*}
By virtue of Vaughan's identity (see \cite[Proposition 13.4]{IK04} for instance), we may decompose the sum over $n$ to bilinear forms
and consider
\begin{align*}
T(\boldsymbol\alpha,\boldsymbol\beta,\boldsymbol\gamma)=\sum_{p\sim P}\gamma_p\mathop{\sum\sum}_{\substack{m\sim M,n\sim N\\(mn,p)=1}}\alpha_m\beta_n\mu^2(mn)\nu(mn)\sym_k(\cos\theta_p(\overline{(mn)^2})),\end{align*}
where $\boldsymbol\alpha=(\alpha_m),\boldsymbol\beta=(\beta_n)$ are some coefficients 
supported in $]M,2M]$ and $]N,2N],$ respectively, 
such that
$|\alpha_m\beta_n|\leqslant10+\log m\log n$. Here $M,N$ are chosen subject to
\begin{align}\label{eq:MN-sizes}
X\cL^{-C}<MN\leqslant X,\ \ \ M\geqslant N,
\end{align}
where $C$ is some large constant. We would like to prove that
\begin{align}\label{eq:T(alpha,beta,gamma)-estimate}
T(\boldsymbol\alpha,\boldsymbol\beta,\boldsymbol\gamma)
&\ll k^A\cL^{A-2}\{PX^{\frac{1}{2}}+P^{\frac{1}{4}}X+P^{\frac{1}{2}}X\cL^{-2A}+(PX)^{\frac{3}{4}}\}\|\boldsymbol\gamma\|,\end{align}
subject to the restrictions in \eqref{eq:MN-sizes}, for any $A>\kappa+2$ and some $C>0$. The lemma then follows from \eqref{eq:T(alpha,beta,gamma)-estimate} immediately.

The restriction $MN>X\cL^{-C}$ is reasonable, since the contributions from those $MN\leqslant X\cL^{-C}$ contribute at most $O(\|\boldsymbol\gamma\|X(P\cL^{\kappa-C})^{\frac{1}{2}})$.
The restriction $M\geqslant N$ is input due to the symmetric roles between $\boldsymbol\alpha$ and $\boldsymbol\beta$.
There is an implicit restriction that $(m,n)=1$ in the inner sums due to the appearance $\mu^2(mn)$, in which case we have $\nu(mn)=\nu(m)\nu(n).$ In this way, we may write
\begin{align*}
T(\boldsymbol\alpha,\boldsymbol\beta,\boldsymbol\gamma)=\sum_{p\sim P}\gamma_p\mathop{\sum\sum}_{\substack{m\sim M,n\sim N\\(mn,p)=(m,n)=1}}\alpha^*(m)\beta^*(n)\sym_k(\cos\theta_p(\overline{(mn)^2})),\end{align*}
with $\alpha^*(m)=\mu^2(m)\alpha_m\nu(m)$ and $\beta^*(n)=\mu^2(n)\beta_n\nu(n)$.
Furthermore, the M\"obius formula gives
\begin{align}\label{eq:trilinear-bilinear}
T(\boldsymbol\alpha,\boldsymbol\beta,\boldsymbol\gamma)=\sum_d\mu(d)\sum_{\substack{p\sim P\\ p\nmid d}}\gamma_p\mathop{\sum\sum}_{\substack{m\sim M/d,n\sim N/d\\(mn,p)=1}}\alpha^*(md)\beta^*(nd)\sym_k(\cos\theta_p(\overline{d^4(mn)^2})).\end{align}
For each fixed $d$, we have two alternative ways to estimate the trilinear forms in \eqref{eq:trilinear-bilinear} by appealing to Lemmas \ref{lm:bilinear} and \ref{lm:bilinearform-avergeovermoduli}.

If $N\leqslant \cL^C$, we then have $M>X\cL^{-2C}$ by
\eqref{eq:MN-sizes}, and Lemma \ref{lm:bilinearform-avergeovermoduli} we may derive that
\begin{align}
T(\boldsymbol\alpha,\boldsymbol\beta,\boldsymbol\gamma)
&\ll k^B\|\boldsymbol\gamma\|\sum_d((M/d)^{\frac{1}{2}}+P)\Big(\sum_{m\sim M/d}|\alpha^*(md)|^2\Big)^{\frac{1}{2}}\Big(\sum_{n\sim N/d}|\beta^*(nd)|\Big)\nonumber\\
&\ll k^B(M^{\frac{1}{2}}+P)M^{\frac{1}{2}}N\cL^{1+\kappa}\|\boldsymbol\gamma\|\nonumber\\
&\ll k^B(X+PX^{\frac{1}{2}}\cL^C)\cL^{1+\kappa}\|\boldsymbol\gamma\|.\label{eq:N-small}\end{align}
We now consider the case $N> \cL^C$. By
\eqref{eq:MN-sizes}, we have $M>X^{\frac{1}{2}}\cL^{-\frac{C}{2}}$. From Lemma \ref{lm:bilinear} it follows that
\begin{align}
T(\boldsymbol\alpha,\boldsymbol\beta,\boldsymbol\gamma)
&\ll k^B(MN)^{\frac{1}{2}}\sum_d\frac{1}{d}\sum_{p\sim P}|\gamma_p|\Big(\sum_{m\sim M/d}|\alpha^*(md)|^2\Big)^{\frac{1}{2}}\Big(\sum_{n\sim N/d}|\beta^*(nd)|^2\Big)^{\frac{1}{2}}\nonumber\\
&\ \ \ \ \times(p^{-\frac{1}{4}}+(d/N)^{\frac{1}{2}}+(d/M)^{\frac{1}{2}}p^{\frac{1}{4}}(\log p)^{\frac{1}{2}})\nonumber\\
&\ll k^BP^{\frac{1}{2}}X\cL^{\kappa+3}(P^{-\frac{1}{4}}+N^{-\frac{1}{2}}+M^{-\frac{1}{2}}P^{\frac{1}{4}})\|\boldsymbol\gamma\|\nonumber\\
&\ll k^BP^{\frac{1}{2}}X\cL^{\kappa+3}(P^{-\frac{1}{4}}+\cL^{-\frac{B}{2}}+(X^{-1}P\cL^C)^{\frac{1}{4}})\|\boldsymbol\gamma\|.\label{eq:N-large}
\end{align}

Combining \eqref{eq:N-small} and \eqref{eq:N-large}, we conclude that
\begin{align*}
T(\boldsymbol\alpha,\boldsymbol\beta,\boldsymbol\gamma)
&\ll k^B\cL^{\kappa+3}\{PX^{\frac{1}{2}}\cL^{C+1}+P^{\frac{1}{4}}X+P^{\frac{1}{2}}X\cL^{-\frac{C}{2}}+(PX)^{\frac{3}{4}}\cL^{\frac{C}{4}}\}\|\boldsymbol\gamma\|\end{align*}
holds uniformly in all tuples $(M,N)$ subject to the restrictions in \eqref{eq:MN-sizes}.
This completes the proof of \eqref{eq:T(alpha,beta,gamma)-estimate}, and thus that of the lemma, by supplying the initial error $O(\|\boldsymbol\gamma\|X(P\cL^{\kappa-C})^{\frac{1}{2}})$ and choosing $A=(C+\kappa+3)/10.$
\endproof

\smallskip

\section{A generalization of the Barban--Davenport--Halberstam theorem}

Regarding the equidistributions of primes in arithmetic progressions, the classical Barban--Davenport--Halberstam theorem (see e.g., \cite[Theorem 17.2]{IK04}) asserts that
\begin{align*}
\sum_{q\leqslant Q}~\sideset{}{^*}\sum_{a\bmod q}\Big|\sum_{\substack{n\leqslant X\\ n\equiv a\bmod q}}\Lambda(n)-\frac{1}{\varphi(q)}\sum_{\substack{n\leqslant X\\ (n,q)=1}}\Lambda(n)\Big|^2\ll X\cL^{-A}\end{align*}
for any $A>0$, as long as $Q\leqslant X\cL^{-B}$ with some $B=B(A)>0,$ where the implied constant depends only on $A$.
As shown by Bombieri, Friedlander and Iwaniec \cite[Theorem 0]{BFI86}, the above estimate also holds if $\Lambda$ is replaced by an arbitrary function $\vartheta_n$ satisfying the  following ``Siegel--Walfisz" condition.

\begin{definition} 
An arithmetic function $\vartheta$ is said to satisfy the ``Siegel--Walfisz" condition, if for any $w\geqslant1,d\geqslant1,(w,a)=1,a\neq0,$
\begin{align}\label{eq:SiegelWalfisz}
\sum_{\substack{n\leqslant X\\ n\equiv a\bmod w\\(n,d)=1}}\vartheta_n-\frac{1}{\varphi(w)}\sum_{\substack{n\leqslant X\\(n,dw)=1}}\vartheta_n\ll\|\boldsymbol\vartheta\|X^{\frac{1}{2}}\tau(d)^{B}\cL^{-A} \end{align}
holds for some constant $B>0$ and any $A>0$ with an implied constant in $\ll$ depending only on $A$.
\end{definition}
In the following treatment to $H_2(X),$ we would require a further generalization, which involves the equidistributions of the convolution of two arbitrary arithmetic functions, and one of them satisfies the 
``Siegel--Walfisz" condition. Moreover, we also require the following definition of admissibility, which concerns with the $q$-analogue of the Mellin transform of $W_q:\bZ/q\bZ\rightarrow\bC,$ defined by
\begin{align*}
\widetilde{W_q}(\chi)
=\frac{1}{\sqrt{q}}~\sideset{}{^*}\sum_{r\bmod q}\overline{\chi}(r)W_q(r).
\end{align*}
Here $\chi$ is a Dirichlet character mod $q.$

\begin{definition} \label{def:admissibility}
Let $q\geqslant1$ a squarefree number, $k\in\bZ$ and $C>0$ a constant. An arithmetic function 
$\varXi_q:\bZ/q\bZ\rightarrow\bC$ is said to be $(k,C)$-admissible, if
\begin{itemize}
\item $\varXi_{q_1q_2}(\cdot)=\varXi_{q_1}(q_2^k\cdot)\varXi_{q_2}(q_1^k\cdot)$ for all $q_1,q_2\geqslant1$ with $\mu^2(q_1q_2)=1;$
\item for each primitive character $\chi\bmod q$, one has $\|\varXi_q\|_\infty+|\widetilde{\varXi}_q(\chi)|\leqslant (\tau(q)\log 2q)^C.$
\end{itemize}
\end{definition}

\begin{remark}
By Lemma \ref{lm:Kloosterman-Mellin-composite}, one may see $\varXi_q$ is $(1,B)$-admissible for some $B>0$
if taking
\begin{align*}
\varXi_q(a)=
\begin{cases}
|\varOmega(a,q)|,\ \ &(a,q)=1,\\
0,&(a,q)>1.
\end{cases}
\end{align*}
\end{remark}

For a $(k,C)$-admissible arithmetic function $\varXi_q$ as above, the Chinese remainder theorem yields
\begin{align}\label{eq:Xi-Mellin-multiplicativity}
\widetilde{\varXi}_q(\chi)=\chi_1(q_2)^k\chi_2(q_1)^k\widetilde{\varXi}_{q_1}(\chi_1)\widetilde{\varXi}_{q_2}(\chi_2)
\end{align}
for all $q_1q_2=q,\chi_1\chi_2=\chi$ with $\chi\bmod{q_1}$ and $\chi_2\bmod{q_2}.$

We are now ready to state our generalization of the Barban--Davenport--Halberstam theorem.

\begin{lemma}\label{lm:bilinear-BDH}
Let $M,N,C>0$ and $q\geqslant1$ squarefree. Let $\boldsymbol\alpha=(\alpha_m)$ be a complex coefficient with support in $[M,2M]$ and also satisfy the  above``Siegel--Walfisz" condition, and $\boldsymbol\beta=(\beta_n),\boldsymbol\gamma_q=(\gamma_{n,q})$ 
complex coefficients with supports in $[N,2N]$ with $\|\boldsymbol\gamma_q\|_\infty \leqslant (\tau(q)\log2q)^C.$
For a $(k,C)$-admissible arithmetic function $\varXi_q$ with some $k\in\bZ$, put
\begin{align*}
\cE(\boldsymbol\alpha,\boldsymbol\beta,\boldsymbol\gamma_q;q,\varXi_q)&=\mathop{\sum\sum}_{(mn,q)=1}\alpha_m\beta_n\gamma_{n,q}
\varXi_q(mn)-\frac{1}{\varphi(q)}\sideset{}{^*}\sum_{r\bmod q}\varXi_q(r)\mathop{\sum\sum}_{(mn,q)=1}\alpha_m\beta_n\gamma_{n,q}\\
&=\sideset{}{^*}\sum_{r\bmod q}\varXi_q(r)\Big(\mathop{\sum\sum}_{mn\equiv r\bmod q}\alpha_m\beta_n\gamma_{n,q}-\frac{1}{\varphi(q)}\mathop{\sum\sum}_{(mn,q)=1}\alpha_m\beta_n\gamma_{n,q}\Big).\end{align*}

Let $r\geqslant1$ and $M\geqslant N$. For any $A>0$, there exists some constant $B=B(A,C)>0,$ such that
\begin{align*}
\sum_{q\leqslant Q}\mu^2(q)\tau(q)^r|\cE(\boldsymbol\alpha,\boldsymbol\beta,\boldsymbol\gamma_{q};q,\varXi_q)|\ll \|\boldsymbol\alpha\| \|\boldsymbol\beta\|Q(MN)^{\frac{1}{2}}(\log MN)^{-A}\end{align*}
for $Q\leqslant MN(\log MN)^{-B},$
where the implied constant depends only on $A,C$ and $r$.
\end{lemma}

\proof
In what follows, we assume $B_0,B_1,B_2,\cdots,B_{11}$ are some positive constants that we will not specialize their values. Moreover, we always keep $q$ to be squarefree.

By virtue of orthogonality of multiplicative characters, we may write
\begin{align*}
\cE(\boldsymbol\alpha,\boldsymbol\beta,\boldsymbol\gamma_{q};q,\varXi_q)
&=\frac{\sqrt{q}}{\varphi(q)}\sum_{\substack{\chi\bmod{q}\\ \chi\neq\chi_0}}\widetilde{\varXi}_q(\chi)\Big(\sum_m\alpha_m\chi(m)\Big)\Big(\sum_n\beta_n\gamma_{n,q}\chi(n)\Big).
\end{align*}
Each non-trivial character $\chi\bmod q$ is induced by some primitive character $\chi^*\bmod{q^*}$ with $q^*\mid q.$ Since $q$ is squarefree, we then have $(q^*,q/q^*)=1$ automatically. Therefore, by \eqref{eq:Xi-Mellin-multiplicativity}, we obtain
\begin{align*}
\cE(\boldsymbol\alpha,\boldsymbol\beta,\boldsymbol\gamma_{q};q,\varXi_q)
&=\frac{\sqrt{q}}{\varphi(q)}\sum_{q^*q_0=q}~\sideset{}{^*}\sum_{\chi\bmod{q^*}}\chi_0(q^*)^k\chi^*(q_0)^k\widetilde{\varXi}_{q^*}(\chi)\widetilde{\varXi}_{q_0}(\chi_0)\\
&\ \ \ \ \ \times\Big(\sum_{(m,q_0)=1}\alpha_m\chi(m)\Big)
\Big(\sum_{(n,q_0)=1}\beta_n\gamma_{n,q}\chi(n)\Big),
\end{align*}
where $\chi_0$ denotes the trivial character mod $q_0$. 
By Definition \ref{def:admissibility}, we have $\widetilde{\varXi}_{q^*}(\chi)\leqslant (\tau(q^*)\log2q^*)^{B_0}$ and $|\widetilde{\varXi}_{q^0}(\chi_0)|\leqslant \sqrt{q_0}(\tau(q_0)\log2q_0)^{B_0}$. It then follows that
\begin{align}
\sum_{q\leqslant Q}\mu^2(q)\tau(q)^r|\cE(\boldsymbol\alpha,\boldsymbol\beta,\boldsymbol\gamma_{q};q,\varXi_q)|
&\ll
Q^{\frac{1}{2}}(\log Q)^{B_1}\sum_{q_0\leqslant Q}\frac{\sqrt{q_0}}{\varphi(q_0)}\sum_{q\leqslant Q/q_0}\frac{\tau(qq_0)^{B_1}}{\varphi(q)}\nonumber\\
&\ \ \ \ \ \times \sideset{}{^*}\sum_{\chi\bmod q}\Big|\sum_{(m,q_0)=1}\alpha_m\chi(m)\Big|
\Big|\sum_{(n,q_0)=1}\beta_n\gamma_{n,qq_0}\chi(n)\Big|\nonumber\\
&=Q^{\frac{1}{2}}(\log Q)^{B_1}\cdot (S_1+S_2),\label{eq:bilinear-S1+S2}
\end{align}
where $S_1$ and $S_2$ denote the corresponding contributions from $q_0\leqslant Q_1$ and $Q_1<q_0\leqslant Q,$ respectively.

By Cauchy's inequality, we find
\begin{align*}
S_1^2\leqslant S_{11}S_{12}\end{align*}
with
\begin{align*}
S_{11}&=\sum_{q_0\leqslant Q_1}\frac{1}{\varphi(q_0)}\sum_{q\leqslant Q/q_0}\frac{1}{\varphi(q)}\sideset{}{^*}\sum_{\chi\bmod q}\Big|\sum_{(m,q_0)=1}\alpha_m\chi(m)\Big|^2,\\
S_{12}&=\sum_{q_0\leqslant Q_1}\frac{q_0}{\varphi(q_0)}\sum_{q\leqslant Q/q_0}\frac{\tau(qq_0)^{2B_1}}{\varphi(q)}\sideset{}{^*}\sum_{\chi\bmod q}\Big|\sum_{(n,q_0)=1}\beta_n\gamma_{n,qq_0}\chi(n)\Big|^2.\end{align*}

We first consider $S_{11}$. We further split $S_{11}$ according to $q\leqslant Q_2$ and $q>Q_2$, and the corresponding contributions are denoted by $S_{11}'$ and $S_{11}'',$ respectively. Regarding $S_{11}'$, the Siegel--Walfisz condition for $\boldsymbol\alpha$ gives
\begin{align*}
S_{11}'&\ll\sum_{q_0\leqslant Q_1}\frac{\tau(q_0)^{B}}{\varphi(q_0)}\sum_{q\leqslant Q_2}\varphi(q)^2\|\boldsymbol\alpha\|^2M(\log M)^{-A}\ll \|\boldsymbol\alpha\|^2Q_2^3M(\log M)^{-A}(\log Q)^{B_2}.\end{align*}
For $S_{11}'',$ the dyadic device yields
\begin{align*}
S_{11}''&\ll\log Q\sum_{q_0\leqslant Q_1}\frac{1}{\varphi(q_0)}\sup_{Q_2<Q_3\leqslant Q/q_0}\frac{1}{Q_3}\sum_{q\sim Q_3}\frac{q}{\varphi(q)}\sideset{}{^*}\sum_{\chi\bmod q}\Big|\sum_{(m,q_0)=1}\alpha_m\chi(m)\Big|^2.\end{align*}
From the classical multiplicative large sieve inequality (see \cite[Theorem 7.13]{IK04} for instance), it follows that
\begin{align*}
S_{11}''&\ll\|\boldsymbol\alpha\|^2\log Q\sum_{q_0\leqslant Q_1}\frac{1}{\varphi(q_0)}\sup_{Q_2<Q_3\leqslant Q/q_0}\frac{1}{Q_3}(Q_3^2+M)\\
&\ll\|\boldsymbol\alpha\|^2(Q+M/Q_2)(\log Q)^2.\end{align*}
Collecting the above estimates for $S_{11}'$ and $S_{11}'',$ we find
\begin{align*}
S_{11}&\ll \|\boldsymbol\alpha\|^2\{Q_2^3M(\log M)^{-A}(\log Q)^{B_2}+(Q+M/Q_2)(\log Q)^2\}.\end{align*}
Taking $Q_2=(\log M)^{A/6}$, we then obtain
\begin{align*}
S_{11}&\ll \|\boldsymbol\alpha\|^2\{M(\log M)^{-A}+Q(\log Q)^{2}\}\end{align*}
by re-defining $A$.

On the other hand, 
\begin{align*}
S_{12}
&\leqslant\sum_{q_0\leqslant Q_1}\frac{q_0}{\varphi(q_0)}\sum_{q\leqslant Q/q_0}\frac{\tau(qq_0)^{2B_1}}{\varphi(q)}\sum_{\chi\bmod q}\Big|\sum_{(n,q_0)=1}\beta_n\gamma_{n,qq_0}\chi(n)\Big|^2\\
&\ll(\log Q)^{B_3}\sum_{q_0\leqslant Q_1}\frac{q_0\tau(q_0)^{B_3}}{\varphi(q_0)}\sum_{q\leqslant Q/q_0}\tau(q)^{B_3}\mathop{\sum\sum}_{n_1\equiv n_2\bmod q}|\beta_{n_1}\beta_{n_2}|.\end{align*}
Note that
\begin{align*}
\sum_{q\leqslant Q/q_0}\tau(q)^{B_3}\mathop{\sum\sum}_{n_1\equiv n_2\bmod q}|\beta_{n_1}\beta_{n_2}|
&\ll \sum_{q\leqslant Q/q_0}\tau(q)^{B_3}\mathop{\sum\sum}_{n_1\equiv n_2\bmod q}|\beta_{n_1}|^2\\
&\ll \|\boldsymbol\beta\|^2Qq_0^{-1}(\log Q)^{B_4}+\sum_{n_1}|\beta_{n_1}|^2\sum_{\substack{n_2\sim N\\n_2\neq n_1}}\tau(|n_2-n_1|)^{B_4}\\
&\ll \|\boldsymbol\beta\|^2(Q/q_0+N)(\log QN)^{B_5},\end{align*}
from which we conclude that
\begin{align*}
S_{12}
&\leqslant \|\boldsymbol\beta\|^2(Q+Q_1N)(\log QN)^{B_6}.\end{align*}

Combining the above estimates for $S_{11}$ and $S_{12},$ we obtain
\begin{align*}
S_1\ll \|\boldsymbol\alpha\|\|\boldsymbol\beta\|(M(\log M)^{-A}+Q)^{\frac{1}{2}}(Q+Q_1N)^{\frac{1}{2}}(\log QN)^{B_7}.\end{align*}

Again by Cauchy's inequality, we find
\begin{align*}
S_2^2\leqslant S_{21}S_{22}\end{align*}
with
\begin{align*}
S_{21}&=\sum_{Q_1<q_0\leqslant Q}\frac{\sqrt{q_0}}{\varphi(q_0)}\sum_{q\leqslant Q/q_0}\frac{\tau(qq_0)^{2B_1}}{\varphi(q)}\sideset{}{^*}\sum_{\chi\bmod q}\Big|\sum_{(m,q_0)=1}\alpha_m\chi(m)\Big|^2,\\
S_{22}&=\sum_{Q_1<q_0\leqslant Q}\frac{\sqrt{q_0}}{\varphi(q_0)}\sum_{q\leqslant Q/q_0}\frac{1}{\varphi(q)}\sideset{}{^*}\sum_{\chi\bmod q}\Big|\sum_{(n,q_0)=1}\beta_n\gamma_{n,qq_0}\chi(n)\Big|^2.\end{align*}
As argued in estimating $S_{11}$ and $S_{12},$ we may derive that
\begin{align*}
S_{21}&\ll \|\boldsymbol\alpha\|^2\{M\sqrt{Q}(\log M)^{-A}+Q(\log Q)^{B_8}\},\\
S_{22}&\ll \|\boldsymbol\beta\|^2\Big(\frac{Q}{\sqrt{Q_1}}+\sqrt{Q}N\Big)(\log QN)^{B_9}.\end{align*}
Therefore, we arrive at
\begin{align*}
S_2\ll \|\boldsymbol\alpha\|\|\boldsymbol\beta\|(M\sqrt{Q}(\log M)^{-A}+Q)^{\frac{1}{2}}\Big(\frac{Q}{\sqrt{Q_1}}+\sqrt{Q}N\Big)^{\frac{1}{2}}(\log QN)^{B_{10}}.\end{align*}

Inserting the estimates for $S_1,S_2$ into \eqref{eq:bilinear-S1+S2}, we find
\begin{align*}
\sum_{q\leqslant Q}\mu^2(q)\tau(q)^r|\cE(\boldsymbol\alpha,\boldsymbol\beta,\boldsymbol\gamma_q;q,\varXi_q)|&\ll \|\boldsymbol\alpha\|\|\boldsymbol\beta\|(MN)^{\frac{1}{2}}Q(\log MNQ)^{B_{11}}\varDelta(M,N,Q;Q_1),\end{align*}
where
\begin{align*}
\varDelta(M,N,Q;Q_1)^2&=\frac{Q}{MN}+\frac{\sqrt{Q}}{M}+\frac{Q_1}{M}+\Big(1+\frac{1}{N}\sqrt{\frac{Q}{Q_1}}\Big)(\log MN)^{-A}.\end{align*}
Taking
\begin{align*}
Q_1=
\begin{cases}
(M/N)^{\frac{2}{3}}Q^{\frac{1}{3}},\ \ \ & M\leqslant NQ,\\
Q(\log MN)^{-A},&M>NQ,
\end{cases}
\end{align*}
the proof is completed by noting that $\sqrt{Q}\leqslant \sqrt{MN(\log MN)^{-B}}\leqslant M(\log MN)^{-B/2}$.
\endproof

\smallskip

\section{Lower bound for $H_1(X)$}

Recalling the definition \eqref{eq:Omega}, we may write
\begin{align*}H_1(X)&= \sum_{n\geqslant1}\varPsi\Big(\frac{n}{X}\Big)\mu^2(n)|\lambda_f(n)-\eta\cdot \varOmega(1,n)|\bigg(\sum_{d|(n,P(z))}\varrho_d\bigg)^2.\end{align*}
To seek a positive lower bound for $H_1(X),$
we
need only consider
those $n$ with few prime factors. To that end, we introduce the interval
\begin{align*}I(P)=~]P,P+P\mathcal{L}^{-1}],\end{align*}
and the set of the products of primes
\begin{align*}\mathcal{P}_i(X;P_{i1},P_{i2},\cdots,P_{ii})&=\{p_1p_2\cdots p_i:p_j\in I(P_{ij})\text{ for each } j\leqslant i\}\end{align*}
for each positive integer $i\geqslant2.$ Furthermore, for each
fixed $i$, we assume that $\{P_{ij}\}$ is a decreasing sequence as powers of $(1+\cL^{-1})$
as $j$ varies, and the product of $P_{ij}$'s falls into $[X,2X]$; i.e.,
\begin{align}\label{eq:Pij-initial}
P_{ij}\exp(-\sqrt{\mathcal{L}})>P_{i(j+1)}\geqslant X^{\frac{1}{12}}\ \ \ (1\leqslant j<i),\ \ \prod_{1\leqslant j\leqslant i}P_{ij}\in[X,2X].\end{align}
In this way, we can bound $H_1(X)$ from below by the summation
over $\mathcal{P}_i(X;P_{i1},P_{i2},\cdots,P_{ii})$; for this, we employ the variants of 
the Sato--Tate distributions stated above.
Due to the positivity of each term, we can drop those $n$'s with
``bad''
arithmetic structures. To this end, we introduce the
following
restrictions on the size of $P_{ij}$:
\begin{equation}\label{eq:Pij-restrictions}
\begin{split}\begin{cases}P_{21}^{\frac{3}{4}}X^\delta<P_{22},\ \ \delta=10^{-2018},\\
\sqrt{P_{31}}\exp(\sqrt{\mathcal{L}})<P_{32},\\
\sqrt{P_{41}}\exp(\sqrt{\mathcal{L}})<P_{42}P_{43},\\
\sqrt{P_{i1}}\exp(\sqrt{\mathcal{L}})<P_{i2}\cdots P_{i(i-1)}\text{ and }\sqrt{P_{i3}\cdots P_{ii}}\exp(\sqrt{\mathcal{L}})<P_{i2},\ \ i\geqslant5.
\end{cases}\end{split}\end{equation}
Now summing up to $i=7,$ we then have the lower bound
\begin{align}\label{eq:H1(X)-lowerbound-initial}
H_1(X)&\geqslant\sum_{2\leqslant i\leqslant7}H_{1,i}(X),
\end{align}
where
\begin{align*}H_{1,i}(X)&=\sideset{}{^\dagger}\sum_{P_{i1},P_{i2},\cdots,P_{ii}}\ \ \sum_{n\in\mathcal{P}_i(X;P_{i1},P_{i2},\cdots,P_{ii})}\varPsi\Big(\frac{n}{X}\Big)|\lambda_f(n)-\eta\cdot\varOmega(1,n)|\bigg(\sum_{d|(n,P(z))}\varrho_d\bigg)^2,\end{align*}
with $\dagger$ yielding that $P_{ij}$'s are powers of $(1+\cL^{-1})$ satisfying the restrictions \eqref{eq:Pij-initial} and \eqref{eq:Pij-restrictions}.

Recalling the choice (\ref{eq:varrho_d}), we find, for each $n\in\mathcal{P}_i(X;P_{i1},P_{i2},\cdots,P_{ii})$, that $n$ has no prime factors smaller than $X^{\frac{1}{12}}$. Note that $z$ will be chosen such that $z\leqslant X^{\frac{1}{12}}$, we then have $(n,P(z))=1$ and
\begin{align*}\sum_{d|(n,P(z))}\varrho_d=\varrho_1=1\end{align*}
for such $n$. Hence we can write
\begin{align}\label{eq:H1i(X)-lowerbound}
H_{1,i}(X)&=(1+o(1))\mathcal{L}^{2i-1}\int_{\mathcal{R}_i}\varSigma(X,\boldsymbol\alpha_i)\ud\boldsymbol\alpha_i,\end{align}
where for $\boldsymbol\alpha_i=(\alpha_2,\cdots,\alpha_i),$ we adopt the convention
\begin{align*}
\cP_i(X,\boldsymbol\alpha_i)
&=\cP_i(X;X^{1-\alpha_2-\cdots-\alpha_i},X^{\alpha_2},\cdots,X^{\alpha_i}),\\
\varSigma(X,\boldsymbol\alpha_i)
&=\sum_{n\in\cP_i(X,\boldsymbol\alpha_i)}\varPsi\Big(\frac{n}{X}\Big)|\lambda_f(n)-\eta\cdot\varOmega(1,n)|,
\end{align*}
and
the multiple-integral is over the area $\mathcal{R}_i$ with
\begin{equation}\label{eq:Ri}
\begin{split}
\mathcal{R}_2&:=\{\alpha_2\in[\tfrac{1}{12},1[:\tfrac{3}{4}(1-\alpha_2)+\delta<\alpha_2<\tfrac{1}{2}\},\\
\mathcal{R}_3&:=\{(\alpha_2,\alpha_3)\in[\tfrac{1}{12},1[^2:\tfrac{1}{2}(1-\alpha_2-\alpha_3)<\alpha_2,\alpha_3<\alpha_2<1-\alpha_2-\alpha_3\},\\
\mathcal{R}_4&:=\{(\alpha_2,\alpha_3,\alpha_4)\in[\tfrac{1}{12},1[^3:\tfrac{1}{2}(1-\alpha_2-\alpha_3-\alpha_4)<\alpha_2+\alpha_3\}\\
&\ \ \ \ \ \cap\{(\alpha_2,\alpha_3,\alpha_4)\in[\tfrac{1}{12},1[^3:\alpha_4<\alpha_3<\alpha_2<1-\alpha_2-\alpha_3-\alpha_4\},\\
\mathcal{R}_i&:=\{(\alpha_2,\cdots,\alpha_i)\in[\tfrac{1}{12},1[^{i-1}:\tfrac{1}{2}(1-\alpha_2-\cdots-\alpha_i)<\alpha_2+\cdots+\alpha_{i-1}\}\\
&\ \ \ \ \ \cap\{(\alpha_2,\cdots,\alpha_i)\in[\tfrac{1}{12},1[^{i-1}:\tfrac{1}{2}(\alpha_3+\cdots+\alpha_i)<\alpha_2\}\\
&\ \ \ \ \ \cap\{(\alpha_2,\cdots,\alpha_i)\in[\tfrac{1}{12},1[^{i-1}:\alpha_i<\alpha_{i-1}<\cdots<\alpha_2<1-\alpha_2-\cdots-\alpha_i\}
\end{split}\end{equation}
for $i\geqslant5$ with $\delta=10^{-2018}$. Note that $\alpha_j<1/j$ for $2\leqslant j\leqslant i$ in the above coordinates.

It remains to seek a lower bound for $\varSigma(X,\boldsymbol\alpha_i)$.
\begin{proposition}\label{prop:Sigma(X,alphai)-lowerbound}
For $i\in[2,7]\cap\bZ$ and $\boldsymbol\alpha_i:=(\alpha_2,\cdots,\alpha_i)\in\mathcal{R}_i,$ we have
\begin{align*}\varSigma(X,\boldsymbol\alpha_i)\geqslant\sqrt{l_i^3/{u_i}}\cdot (1+o(1))\cdot |\cP_i(X,\boldsymbol\alpha_i)|\end{align*}
for all sufficiently large $X,$ where 
\begin{align*}l_i=(1-4|\eta|\cdot(\tfrac{8}{3\pi})^{i-1}(\tfrac{11}{12})^i+B_i\eta^2)^+\end{align*}
and
\begin{align*}
u_i=16^i\cdot|\eta|^4+4\cdot (\tfrac{22}{3})^i\cdot |\eta|^3+6\cdot 4^i\cdot |\eta|^2+4\cdot (2\sqrt{5})^i\cdot |\eta|+2^i\end{align*}
with the convention that $x^+=\max\{0,x\}$ and $B_i$'s being given below by $\eqref{eq:Bi-constants}.$

Consequently, for $i\in[2,7]\cap\bZ,$ we have
\begin{align}\label{eq:Sigma(X,alphai)-integrallowerbound}
\int_{\mathcal{R}_i}\varSigma(X,\boldsymbol\alpha_i)\ud\boldsymbol\alpha_i
&\geqslant (1+o(1))\sqrt{l_i^3/{u_i}}\cdot I_i,\end{align}
where
\begin{align}\label{eq:Ii}
I_i=\int_{\mathcal{R}_i}\frac{\ud\boldsymbol\alpha_i}{\alpha_2\cdots \alpha_i(1-\alpha_2-\cdots-\alpha_i)}.
\end{align}
\end{proposition}

The proof of Proposition \ref{prop:Sigma(X,alphai)-lowerbound} will be given in the next section. Proposition \ref{prop:H1(X)-lowerbound} then follows by substituting \eqref{eq:H1i(X)-lowerbound} and \eqref{eq:Sigma(X,alphai)-integrallowerbound} into \eqref{eq:H1(X)-lowerbound-initial}.

\smallskip

\section{Proof of Proposition \ref{prop:Sigma(X,alphai)-lowerbound}}

For the seek of proving Proposition \ref{prop:Sigma(X,alphai)-lowerbound}, we would like to introduce the following averages
\begin{align*}
\cA_{\ell}(X,\boldsymbol\alpha_i)&=\sum_{n\in\cP_i(X,\boldsymbol\alpha_i)}\varPsi\Big(\frac{n}{X}\Big)|\lambda_f(n)|^\ell,\\
\cB(X,\boldsymbol\alpha_i)&=\sum_{n\in\cP_i(X,\boldsymbol\alpha_i)}\varPsi\Big(\frac{n}{X}\Big)|\varOmega(1,n)|^2,\\
\cC(X,\boldsymbol\alpha_i)&=\sum_{n\in\cP_i(X,\boldsymbol\alpha_i)}\varPsi\Big(\frac{n}{X}\Big)\lambda_f(n)\varOmega(1,n),\end{align*}
where $\cP_i(X,\boldsymbol\alpha_i)$ is given as before for $i\in[2,7]\cap\bZ.$

By H\"older's inequality, we have
\begin{align*}
\varSigma(X,\boldsymbol\alpha_i)\geqslant \varSigma_{2}(X,\boldsymbol\alpha_i)^{\frac{3}{2}}\varSigma_{4}(X,\boldsymbol\alpha_i)^{-\frac{1}{2}},\end{align*}
where
\begin{align*}\varSigma_{\ell}(X,\boldsymbol\alpha_i):=\sum_{n\in\cP_i(X,\boldsymbol\alpha_i)}\varPsi\Big(\frac{n}{X}\Big)|\lambda_f(n)-\eta\cdot\varOmega(1,n)|^\ell.\end{align*}

To prove Proposition \ref{prop:Sigma(X,alphai)-lowerbound}, it suffices to prove that
\begin{align*}\varSigma_{2}(X,\boldsymbol\alpha_i)\geqslant l_i(1+o(1))\cdot |\cP_i(X,\boldsymbol\alpha_i)|,\end{align*}
\begin{align*}\varSigma_{4}(X,\boldsymbol\alpha_i)\leqslant u_i(1+o(1))\cdot |\cP_i(X,\boldsymbol\alpha_i)|\end{align*}
with $l_i,u_i$ are given as in Proposition \ref{prop:Sigma(X,alphai)-lowerbound}.

\subsection{Bounding $\varSigma_{4}(X,\boldsymbol\alpha_i)$ from above}
From Weil's bound for Kloosterman sums, we have
\begin{align}\label{eq:varSigma4(X,alphai)-upperbound-initial}
\varSigma_{4}(X,\boldsymbol\alpha_i)
&\leqslant \sum_{0\leqslant \ell\leqslant4}\binom{4}{\ell}\cdot(2^i\cdot|\eta|)^{4-\ell}\cdot \cA_{\ell}(X,\boldsymbol\alpha_i).
\end{align}
By the definition of $\cP_i(X,\boldsymbol\alpha_i)$ and multiplicativity of Hecke eigenvalues, it follows from Lemma \ref{lm:pi-kappa(X)-moments} that
\begin{align*}\cA_{\ell}(X,\boldsymbol\alpha_i)\leqslant c_\ell^{i}(1+o(1))|\cP_i(X,\boldsymbol\alpha_i)|,\end{align*}
from which and \eqref{eq:varSigma4(X,alphai)-upperbound-initial} we conclude that
\begin{align*}\varSigma_{4}(X,\boldsymbol\alpha_i)
&\leqslant \sum_{0\leqslant \ell\leqslant4}\binom{4}{\ell}\cdot(2^i\cdot|\eta|)^{4-\ell}\cdot c_\ell^{i}\cdot (1+o(1))|\cP_i(X,\boldsymbol\alpha_i)|\\
&=u_i\cdot (1+o(1)) |\cP_i(X,\boldsymbol\alpha_i)|,\end{align*}
provided that $X$ is large enough, where $u_i$'s are given as claimed.

\subsection{Bounding $\varSigma_{2}(X,\boldsymbol\alpha_i)$ from below}
We now turn to the lower bound for $\varSigma_{2}(X,\boldsymbol\alpha_i).$
Squaring out, we may write
\begin{align}\label{eq:Sigma2(X,alphai)}
\varSigma_{2}(X,\boldsymbol\alpha_i)
&=\cA_2(X,\boldsymbol\alpha_i)+\eta^2\cdot\cB(X,\boldsymbol\alpha_i)-2\eta\cdot\cC(X,\boldsymbol\alpha_i).\end{align}

By the definition of $\cP_i(X,\boldsymbol\alpha_i)$ and multiplicativity of Hecke eigenvalues, if follows from Lemma \ref{lm:pi-kappa(X)-moments} that
\begin{align}\label{eq:A2(X,alphai)}
\cA_2(X,\boldsymbol\alpha_i)=(1+o(1))|\cP_i(X,\boldsymbol\alpha_i)|.
\end{align}

The lower bound for $\cB(X,\boldsymbol\alpha_i)$
follows from  {\it joint} equidistributions of Kloosterman sums. By twisted multiplicatitivity, $\varOmega(1,n)$ can be expressed as the product the two Kloosterman sums, the  equidistributions of which are known in a certain sense. To formulate the precise distributions, we would like to introduce the corresponding measures firstly. Following \cite{FM07}, we define a measure $\mu^{(1)}$ on $[-1,1]$ to be the image of the measure $\mu_{\mathrm{ST}}$ under the mapping $\theta\mapsto \cos\theta,$ so that $\ud\mu^{(1)}=\frac{2}{\pi}\sqrt{1-x^2}\ud x$. 
Furthermore, for $k\geqslant2$, we define a measure $\mu^{(k)}$ on $[-1,1]$ to be the image of $\mu^{(1)}\otimes \mu^{(1)}\otimes\cdots\otimes \mu^{(1)}$ under the mapping 
\begin{align*}
[-1,1]^k&\rightarrow [-1,1]\\
(x_1,x_2,\cdots,x_k)&\mapsto x_1x_2\cdots x_k.
\end{align*}
Then for $x\in[0,1],$ we have the following recursive relation
\begin{align}\label{eq:mu.1}
\mu^{(1)}([-x,x])=\frac{4}{\pi}\int_0^x\sqrt{1-t^2}\ud t,
\end{align}
\begin{align}\label{eq:mu.k}
\mu^{(k)}([-x,x])=\mu^{(1)}([-x,x])+\frac{4}{\pi}\int_x^1\mu^{(k-1)}([-x/t,x/t])\sqrt{1-t^2}\ud t,\ \ \ k\geqslant2.
\end{align}

\begin{lemma}\label{lm:equidistribution}
With the notation as above, for $i\in[2,7]\cap\bZ$ and $\boldsymbol\alpha_i:=(\alpha_2,\cdots,\alpha_i)\in\mathcal{R}_i$ as given by $\eqref{eq:Ri},$ the sets
\[\{2^{1-i}\varOmega(p_1,p_2\cdots p_i):n=p_1p_2\cdots p_i\in \cP_i(X,\boldsymbol\alpha_i)\}\]
and
\[\{2^{1-i}\varOmega(p_2\cdots p_i,p_1):n=p_1p_2\cdots p_i\in\cP_i(X,\boldsymbol\alpha_i)\}\]
equidistribute in $[-1,1]$ with respect to $\mu^{(i-1)}$ and
$\mu^{(1)}$,
respectively, as $X\rightarrow+\infty$, where the measures $\mu^{(j)}$ on $[-1,1]$ are defined recursively by
$\eqref{eq:mu.1}$ and $\eqref{eq:mu.k}.$
\end{lemma}

The original statement of Lemma \ref{lm:equidistribution}, in the case $i\in\{3,4,5\}$, can be found in \cite[Propositions 6.1, 6.2 and 6.3]{FM07} and the case $i\in\{6,7\}$ can be treated in a similar way.
The case $i=2$ follows from \cite[Theorem 1.5]{FKM14}
by taking $K(n)=\sym_k(\theta_p(\overline{n^2}))$ therein.

The following rearrangement type inequality, due to Matom\"{a}ki \cite{Ma11}, allows us to derive a lower bound for $\cB(X,\boldsymbol\alpha_i)$ from the equidistributions of Kloosterman sums arising from the above factorization.
\begin{lemma}\label{lm:Matomaki}
Assume that the sequences $(a_n)_{n\leqslant N}$ and $(b_n)_{n\leqslant N}$ contained
in $[0,1]$  equidistribute with respect to some absolutely continuous measures $\mu_a$ and $\mu_b$, respectively, as $N\rightarrow+\infty$. Then
\begin{align*}(1+o(1))\int_0^1xy_l(x)\ud\mu_a([0,x]) \leqslant\frac{1}{N}\sum_{n\leqslant N}a_nb_n\leqslant
(1+o(1))\int_0^1xy_u(x)\ud\mu_a([0,x]),\end{align*}
where $y_l(x)$ is the smallest solution to the equation $\mu_b([y_l,1])=\mu_a([0,x])$ and
$y_u(x)$ is the largest solution to the equation $\mu_b([0,y_u])=\mu_a([0,x])$.
\end{lemma}

We now write
\begin{align*}
\cB(X,\boldsymbol\alpha_i)
&=\sum_{p_1p_2\cdots p_i\in\cP_i(X,\boldsymbol\alpha_i)}\varPsi\Big(\frac{p_1p_2\cdots p_i}{X}\Big)|\varOmega(p_2\cdots p_i,p_1)|^2|\varOmega(p_1,p_2\cdots p_i)|^2.\end{align*}
By Lemma \ref{lm:Matomaki}, we have
\begin{align}\label{eq:jointequidistribution-lowerbound}
\cB(X,\boldsymbol\alpha_i)\geqslant B_i(1+o(1))|\cP_i(X,\boldsymbol\alpha_i)|\end{align}
with
\begin{align*}B_i=4^i\int_0^1x^2 y_i(x)^2\ud\mu^{(1)}([-x,x]),\end{align*}
where $y_i(x)$ is the unique solution to the equation
\begin{align*}\mu^{(1)}([-x,x])=\mu^{(i-1)}([-1,-y]\cup[y,1])=1-\mu^{(i-1)}([-y,y]).\end{align*}
With the help of Mathematica 10, we can obtain
\begin{equation}\label{eq:Bi-constants}
\begin{split}
B_2\geqslant 0.233838, \ \ \ \ \ \ \ B_5&\geqslant 0.023523\\
B_3\geqslant 0.099779, \ \ \ \ \ \ \ B_6&\geqslant 0.011685\\
B_4\geqslant 0.047473, \ \ \ \ \ \ \ B_7&\geqslant 0.005567.
\end{split}
\end{equation}

To conclude Proposition \ref{prop:Sigma(X,alphai)-lowerbound}, it remains to control $\cC(X,\boldsymbol\alpha_i)$ effectively.
It is highly desired that $\lambda_f(n)$ does not correlate with $\varOmega(1,n)$ as $n$ runs over primes or almost primes. Quantitatively, we expect, as discussed in Section \ref{sec:outline}, that 
\begin{align*}
\cC(X,\boldsymbol\alpha_i)&=o(|\cP_i(X,\boldsymbol\alpha_i)|)\end{align*}
for $\boldsymbol\alpha_i\in\cR_i$ as given by \eqref{eq:Ri} and $X\rightarrow+\infty.$
Unfortunately, this non-correlation is not yet known even as $n$ runs over consecutive integers.
Our success builds on the observation that $|\lambda_f(p)|^2$ is approximately $1$ on average; however, $|\lambda_f(p)|$ and $|\varOmega(n,p)|$ are both smaller than 1 on average in a suitable family, so that one may obtain a relatively small scalar in the upper bound of $\cC(X,\boldsymbol\alpha_i)$, even though the sign changes of $\lambda_f(n)\varOmega(1,n)$ are not taken into account.

Precisely speaking, we are able to bound $\cC(X,\boldsymbol\alpha_i)$ as follows.
\begin{proposition}\label{prop:C(X,alphai)-upperbound}
With the notation as above, we have, for all sufficiently large $X,$ that
\begin{align*}
|\cC(X,\boldsymbol\alpha_i)|
&\leqslant2\cdot\Big(\frac{8}{3\pi}\Big)^{i-1}\Big(\frac{11}{12}\Big)^{i}(1+o(1))|\cP_i(X,\boldsymbol\alpha_i)|\end{align*}
for each $i\in[2,7]\cap\bZ.$
\end{proposition}

The lower bound for  $\varSigma_{2}(X,\boldsymbol\alpha_i)$ in Proposition \ref{prop:Sigma(X,alphai)-lowerbound} then follows by combining \eqref{eq:Sigma2(X,alphai)}, \eqref{eq:A2(X,alphai)}, \eqref{eq:jointequidistribution-lowerbound} and Proposition
\ref{prop:C(X,alphai)-upperbound}, as well as 
\begin{align*}
|\cP_i(X,\boldsymbol\alpha_i)|=\frac{X\mathcal{L}^{-2i}(1+o(1))}{\alpha_2\cdots\alpha_i(1-\alpha_2-\cdots-\alpha_i)}\end{align*}
from the prime number theorem. The complete proof of Proposition $\ref{prop:C(X,alphai)-upperbound}$ will be given in the next section.

\smallskip

\section{Proof of Proposition $\ref{prop:C(X,alphai)-upperbound}$}
By the definition of $\cP_i(X,\boldsymbol\alpha_i)$ and twisted multiplicativity of Kloosterman sums, we may write
\begin{align*}
\cC(X,\boldsymbol\alpha_i)
&=\sum_{p_1p_2\cdots p_i\in\cP_i(X,\boldsymbol\alpha_i)}\varPsi\Big(\frac{p_1p_2\cdots p_i}{X}\Big)\lambda_f(p_1p_2\cdots p_i)\varOmega(p_2\cdots p_i,p_1)\varOmega(p_1,p_2\cdots p_i).\end{align*}
Weil's bound gives
\begin{align*}
|\cC(X,\boldsymbol\alpha_i)|
&\leqslant2\cC^*(X,\boldsymbol\alpha_i),\end{align*}
where
\begin{align*}
\cC^*(X,\boldsymbol\alpha_i)
&=\sum_{p_1p_2\cdots p_i\in\cP_i(X,\boldsymbol\alpha_i)}\varPsi\Big(\frac{p_1p_2\cdots p_i}{X}\Big)|\lambda_f(p_1p_2\cdots p_i)||\varOmega(p_1,p_2\cdots p_i)|.\end{align*}
It suffices to prove that
\begin{align*}
\cC^*(X,\boldsymbol\alpha_i)
&\leqslant\Big(\frac{8}{3\pi}\Big)^{i-1}\Big(\frac{11}{12}\Big)^{i}(1+o(1))|\cP_i(X,\boldsymbol\alpha_i)|\end{align*}
for $i\in[2,7]\cap\bZ.$
We prove these inequalities case by case. The case $i=2$ is a bit different, which essentially relies on Lemma \ref{lm:bilinearaverageoverprimes} and the remaining cases will be concluded by Lemmas \ref{lm:bilinear} and \ref{lm:bilinearform-twisted} amongst other things.

\subsection{Bounding $\cC^*(X,\boldsymbol\alpha_2)$}
We first consider the case $i=2.$
From the twisted multiplicativity for Kloosterman sums and multiplicativity for Hecke eigenvalues, we may write
\begin{align*}
\cC^*(X,\boldsymbol\alpha_2)
&=\sum_{p_1p_2\in\cP_2(X,\boldsymbol\alpha_2)}\varPsi\Big(\frac{p_1p_2}{X}\Big)|\lambda_f(p_1)||\lambda_f(p_2)||\varOmega(p_1,p_2)|.\end{align*}
We then apply Lemma \ref{lm:bilinearaverageoverprimes} with
\[(n,p)=(p_1,p_2),\ \ \nu(n)=|\lambda_f(p_1)|,\ \ \gamma_q=|\lambda_f(p_2)|,\]
getting
\begin{align*}
\cC^*(X,\boldsymbol\alpha_2)&=\frac{8}{3\pi}\sum_{p_1p_2\in\cP_2(X,\boldsymbol\alpha_2)}\varPsi\Big(\frac{p_1p_2}{X}\Big)|\lambda_f(p_1)||\lambda_f(p_2)|\\
&\ \ \ \ +O\Big(\cL^{10}\{X^{\frac{1}{2}+\alpha_2}+X^{1-\frac{1}{4}\alpha_2}+X\cL^{-20}+X^{\frac{3+2\alpha_2}{4}}\}\Big).
\end{align*}
The desired inequality for $i=2$ now follows from Lemma \ref{lm:pi-kappa(X)-moments}
and $\frac{3}{7}<\alpha_2<\frac{1}{2}$ by \eqref{eq:Ri}.

\subsection{Bounding $\cC^*(X,\boldsymbol\alpha_3)$}
We now consider the case $i=3.$ By multiplicativity, we may write
\begin{align*}
\cC^*(X,\boldsymbol\alpha_3)
&=\sum_{p_1p_2p_3\in\cP_3(X,\boldsymbol\alpha_3)}\varPsi\Big(\frac{p_1p_2p_3}{X}\Big)|\lambda_f(p_1)||\lambda_f(p_2)||\lambda_f(p_3)||\varOmega(p_1p_3,p_2)||\varOmega(p_1p_2,p_3)|.\end{align*}
In view of the Chebyshev approximation for $|\cos\theta|$ (see Lemma \ref{lm:cos-Chebyshev}), we consider
\begin{align*}
\cC_k^*(X,\boldsymbol\alpha_3)
&:=\sum_{p_1p_2p_3\in\cP_3(X,\boldsymbol\alpha_3)}\varPsi\Big(\frac{p_1p_2p_3}{X}\Big)|\lambda_f(p_1)||\lambda_f(p_2)||\lambda_f(p_3)||\varOmega(p_1p_3,p_2)|\sym_{k}(\cos\theta_{p_3}(\overline{(p_1p_2)^2})).\end{align*}
Applying Lemma \ref{lm:bilinearform-twisted} with
\begin{align*}
s=1,\ (m,n,q)=(p_1,p_2,p_3),\ (M,N)=(X^{1-\alpha_2-\alpha_3},X^{\alpha_2}),\end{align*}
\begin{align*}
\alpha_m=|\lambda_f(p_1)|,\ \beta_n=|\lambda_f(p_2)|,\ \gamma_{m,n}=|\varOmega(p_1p_3,p_2)|,\end{align*}
we obtain
\begin{align*}
\cC_k^*(X,\boldsymbol\alpha_3)
&\ll(k+1)\sum_{p_1p_2p_3\in\cP_3(X,\boldsymbol\alpha_3)}|\lambda_f(p_3)|(p_3^{-\frac{1}{8}}+X^{-\frac{\alpha_2}{4}}p_3^{\frac{1}{8}}+X^{\frac{2\alpha_2+\alpha_3-1}{2}})\\
&\ll(k+1)\exp(-\sqrt{\cL})|\cP_3(X,\boldsymbol\alpha_3)|
\end{align*}
by Cauchy's inequality and Lemma \ref{lm:pi-kappa(X)-moments}.
Therefore, it follows from Lemma \ref{lm:cos-Chebyshev} that
\begin{align*}
\cC^*(X,\boldsymbol\alpha_3)
&=\frac{8}{3\pi}(1+o(1))\sum_{p_1p_2p_3\in\cP_3(X,\boldsymbol\alpha_3)}\varPsi\Big(\frac{p_1p_2p_3}{X}\Big)|\lambda_f(p_1)||\lambda_f(p_2)||\lambda_f(p_3)||\varOmega(p_1p_3,p_2)|.\end{align*}
By Lemmas \ref{lm:bilinear} and \ref{lm:cos-Chebyshev}, we further have
\begin{align*}
\cC^*(X,\boldsymbol\alpha_3)
&=\Big(\frac{8}{3\pi}\Big)^2(1+o(1))\sum_{p_1p_2p_3\in\cP_3(X,\boldsymbol\alpha_3)}\varPsi\Big(\frac{p_1p_2p_3}{X}\Big)|\lambda_f(p_1)||\lambda_f(p_2)||\lambda_f(p_3)|.\end{align*}
Then Lemma \ref{lm:pi-kappa(X)-moments} yields
\begin{align*}
\cC^*(X,\boldsymbol\alpha_3)
&\leqslant\Big(\frac{8}{3\pi}\Big)^2\Big(\frac{11}{12}\Big)^3(1+o(1))|\cP_3(X,\boldsymbol\alpha_3)|\end{align*}
as expected.

\subsection{Bounding $\cC^*(X,\boldsymbol\alpha_i)$ for $i\in[4,7]\cap\bZ$}
The cases for $i\geqslant4$ can be treated in a similar way to that for $i=3$, and we only present the details for $i=7$ here.
From multiplicativities, we may write
\begin{align*}
\cC^*(X,\boldsymbol\alpha_7)
&=\sum_{p_1p_2\cdots p_7\in\cP_7(X,\boldsymbol\alpha_7)}\varPsi\Big(\frac{p_1p_2\cdots p_7}{X}\Big)|\lambda_f(p_1)||\lambda_f(p_2)|\cdots|\lambda_f(p_7)|\\
&\ \ \ \ \ \times\prod_{2\leqslant j\leqslant7}|\varOmega(p_1p_2\cdots p_7/p_j,p_j)|.\end{align*}
In view of Lemma \ref{lm:cos-Chebyshev}, we consider
\begin{align*}
\cC_{\bk}^*(X,\boldsymbol\alpha_7)
&=\sum_{p_1p_2\cdots p_7\in\cP_7(X,\boldsymbol\alpha_7)}\varPsi\Big(\frac{p_1p_2\cdots p_7}{X}\Big)|\lambda_f(p_1)||\lambda_f(p_2)|\cdots|\lambda_f(p_7)||\varOmega(p_1p_3p_4\cdots p_7,p_2)|\\
&\ \ \ \ \ \ \times\prod_{3\leqslant j\leqslant 7}\sym_{k_{j-2}}(\cos\theta_{p_{j}}(\overline{(p_1p_2\cdots p_7/p_{j})^2}))\end{align*}
for $\bk=(k_1,\cdots,k_5)\in\bZ_{\geqslant0}^5.$
The term with $\bk=(0,\cdots,0)$ is expected to contribute as the main term. We now assume at least one of $k_1,k_2,\cdots,k_5$ is positive, and only consider the case $k_1k_2\cdots k_5\neq0$ without loss of generality (the remaining cases are simpler).
Applying Lemma \ref{lm:bilinearform-twisted} with
\begin{align*}
s=5,\ (m,n,q)=(p_1,p_2,p_3p_4\cdots p_7),\ (M,N)=(X^{1-\alpha_2-\cdots-\alpha_7},X^{\alpha_2}),\end{align*}
\begin{align*}
\alpha_m=|\lambda_f(p_1)|,\ \beta_n=|\lambda_f(p_2)|,\ \gamma_{m,n}=|\varOmega(p_1p_3p_4\cdots p_7,p_2)|,\end{align*}
we get
\begin{align*}
\cC_\bk^*(X,\boldsymbol\alpha_7)
&\ll k_1k_2\cdots k_5\sum_{p_1p_2\cdots p_7\in\cP_7(X,\boldsymbol\alpha_7)}|\lambda_f(p_1)||\lambda_f(p_2)|\cdots|\lambda_f(p_7)|\\
&\ \ \ \ \ \times\{(p_3p_4\cdots p_7)^{-\frac{1}{8}}+X^{-\frac{\alpha_2}{4}}(p_3p_4\cdots p_7)^{\frac{1}{8}}+X^{\frac{2\alpha_2+\alpha_3+\alpha_4+\cdots \alpha_7-1}{2}}\}\\
&\ll k_1k_2\cdots k_5\exp(-\sqrt{\cL})|\cP_7(X,\boldsymbol\alpha_7)|
\end{align*}
by Cauchy's inequality and Lemma \ref{lm:pi-kappa(X)-moments}. Therefore, it follows from Lemmas \ref{lm:cos-Chebyshev} and \ref{lm:bilinear} that
\begin{align*}
\cC^*(X,\boldsymbol\alpha_7)
&=\Big(\frac{8}{3\pi}\Big)^6(1+o(1))\sum_{p_1p_2\cdots p_7\in\cP_7(X,\boldsymbol\alpha_7)}\varPsi\Big(\frac{p_1p_2\cdots p_7}{X}\Big)|\lambda_f(p_1)||\lambda_f(p_2)|\cdots|\lambda_f(p_7)|,\end{align*}
which yields the desired upper bound in view of Lemma \ref{lm:pi-kappa(X)-moments}.

\smallskip

\section{Upper bound for $H_2(X)$}
First, we may write
\begin{align}\label{eq:H2(X)-upperbound-initial}
H_2(X)&\leqslant H_{21}(X)+|\eta|\cdot H_{22}(X)
\end{align}
with
\begin{align*}H_{21}(X)&=\sum_{n\geqslant1}\varPsi\Big(\frac{n}{X}\Big)\mu^2(n)|\lambda_f(n)|\tau_{\varDelta}(n;\alpha,\beta)\Big(\sum_{d|(n,P(z))}\varrho_d\Big)^2,\\
H_{22}(X)&=\sum_{n\geqslant1}\varPsi\Big(\frac{n}{X}\Big)\mu^2(n)|\varOmega(1,n)|\tau_{\varDelta}(n;\alpha,\beta)\Big(\sum_{d|(n,P(z))}\varrho_d\Big)^2.\end{align*}

\subsection{Dimension-reduction in $H_{22}(X)$}
We now transform $H_{22}(X)$ in the flavor of \cite{Xi18}, so that the {\it dimension} 
of sifting $|\varOmega(1,n)|$ in $H_{22}(X)$ can be reduced. It is a pity that there are some slips in the original arguments of \cite{Xi18}, which will be definitely remedied in this section. As one may find from the definition \eqref{eq:tau(n;alpha,beta)}, the restriction $p\mid d\Rightarrow p>(\log n)^A$ in \cite{Xi18} is replaced by $d\leqslant n^\frac{1}{1+\varDelta}$ with $\varDelta>1$. This new restriction is technical, and it will be reflected in the application of Lemma \ref{lm:bilinear-BDH}.

By twisted multiplicativity,  $H_{22}(X)$ becomes
\begin{align*}H_{22}(X)&=\mathop{\sum\sum}_{m^\varDelta\leqslant n}\varPsi\Big(\frac{mn}{X}\Big)\mu^2(mn)|\varOmega(m,n)||\varOmega(n,m)|\alpha^{\omega(m)}\beta^{\omega(n)}\Big(\sum_{d|(mn,P(z))}\varrho_d\Big)^2.\end{align*}
The Weil bound gives
\begin{align*}H_{22}(X)&\leqslant\mathop{\sum\sum}_{m^\varDelta\leqslant n}\varPsi\Big(\frac{mn}{X}\Big)\mu^2(mn)|\varOmega(n,m)|\alpha^{\omega(m)}(2\beta)^{\omega(n)}\Big(\sum_{d|(mn,P(z))}\varrho_d\Big)^2.\end{align*}
Following the arguments on smooth partition of units in \cite{Xi18} (see e.g., \cite{Fo85}), we have
\begin{align}\label{eq:H22(X)-partition}
H_{22}(X)&\leqslant \sum_{(M,N)}H_{22}(X;M,N),\end{align}
where $M,N$ run over powers of $1+\cL^{-B}$ with $B$ appropriately large and
\begin{align*}H_{22}(X;M,N)&=\mathop{\sum\sum}_{m^\varDelta\leqslant nd}U(m)V(n)\varPsi\Big(\frac{mn}{X}\Big)\mu^2(mn)|\varOmega(n,m)|\alpha^{\omega(m)}(2\beta)^{\omega(n)}\Big(\sum_{d|(mn,P(z))}\varrho_d\Big)^2\end{align*}
with $U,V$ being certain smooth functions supported on $]M,M(1+\cL^{-B})]$ and 
$]N,N(1+\cL^{-B})]$, respectively. By symmetry, we may assume that
\begin{align}\label{eq:sizes-MN}
MN\asymp X,\ \ \ M^\varDelta\ll N.
\end{align}
Note that there are at most $O(\cL^{2B+2})$ tuples of $(M,N)$ in summation.

To evaluate $H_{22}(X;M,N)$, we make the transformation by
\begin{align*}H_{22}(X;M,N)&=\sum_{m}U(m)\mu^2(m)\alpha^{\omega(m)}\Xi(m,N),\end{align*}
where
\begin{align}\label{eq:xi}
\xi(n)=\mathop{\sum\sum}_{[d_1,d_2]=n}\varrho_{d_1}\varrho_{d_2}
\end{align}
and
\begin{align*}\Xi(m,N):=\mathop{\sum\sum}_{\substack{(nd,m)=1\\m^\varDelta\leqslant nd\\ d\mid P(z)}}V(nd)\varPsi\Big(\frac{mnd}{X}\Big)|\varOmega(nd,m)|\mu^2(nd)(2\beta)^{\omega(nd)}\sum_{l\mid (m,P(z))}\xi(dl).\end{align*}
Moreover, one can employ Mellin inversion to separate variables $n,d$ subject to the restrictions in $nd>m$, $V(nd)$ and $\varPsi(mnd/X)$.
Due to the appearance of $\mu^2(nd)$, we can also introduce the M\"obius formula to relax the implicit restriction $(n,d)=1.$
Noting that $N\geqslant X^{\frac{\varDelta}{1+\varDelta}}\gg\sqrt{X}$ in view of \eqref{eq:sizes-MN} and $\xi$ is supported on squarefree numbers up to $\sqrt{X}\exp(-2\sqrt{\cL})$ by the choice of $(\varrho_d)$, we are in a good position to apply Lemmas \ref{lm:bilinear-BDH}
and \ref{lm:Kloosterman-Mellin-prime}, getting
\begin{align*}H_{22}(X;M,N)&=\sum_mU(m)\mu^2(m)\Big(\frac{8\alpha}{3\pi}\Big)^{\omega(m)}\Xi^*(m,N)+O(X\cL^{-2B-4}),\end{align*}
where
\begin{align*}\Xi^*(m,N):=\mathop{\sum\sum}_{\substack{(nd,m)=1\\m^\varDelta\leqslant n\\ d\mid P(z)}}V(nd)\varPsi\Big(\frac{mnd}{X}\Big)\mu^2(nd)(2\beta)^{\omega(nd)}\sum_{l\mid (m,P(z))}\xi(dl).\end{align*}
Rearranging all above summations, we may obtain
\begin{align*}H_{22}(X;M,N)&=\mathop{\sum\sum}_{m^\varDelta\leqslant n}U(m)V(n)\varPsi\Big(\frac{mn}{X}\Big)\mu^2(mn)\Big(\frac{8\alpha}{3\pi}\Big)^{\omega(m)}(2\beta)^{\omega(n)}\Big(\sum_{d|(mn,P(z))}\varrho_d\Big)^2\\
&\ \ \ \ +O(X\cL^{-2B-4}).\end{align*}

Taking into account all admissible tuples $(M,N)$, we find
\begin{align*}H_{22}(X)&\leqslant\sum_{n}\varPsi\Big(\frac{n}{X}\Big)\mu^2(n)\tau_\varDelta\Big(n;\frac{8\alpha}{3\pi},2\beta\Big)\Big(\sum_{d|(n,P(z))}\varrho_d\Big)^2+O(X\cL^{-2}).\end{align*}
Taking $\alpha,\beta>0$ such that 
\begin{align}\label{eq:alpha-beta1}
\frac{8\alpha}{3\pi}+2\beta\leqslant2,
\end{align}
so that
\begin{align*}\tau_\varDelta\Big(n;\frac{8\alpha}{3\pi},2\beta\Big)\leqslant2^{\omega(n)}\end{align*}
for all squarefree $n\geqslant1$. Hence the above upper bound for $H_{22}(X)$ becomes
\begin{align}\label{eq:H22(X)-upperbound-dimreduction}
H_{22}(X)&\leqslant\frac{1}{2}\sum_{n}\varPsi\Big(\frac{n}{X}\Big)\mu^2(n)2^{\omega(n)}\Big(\sum_{d|(n,P(z))}\varrho_d\Big)^2+O(X\cL^{-2}).\end{align}

\subsection{Bounding $H_{21}(X)$ initially}
On the other hand, from the trivial inequality $\tau_{\varDelta}(n;\alpha,\beta)\leqslant(\alpha+\beta)^{\omega(n)}$ it follows that
\begin{align*}H_{21}(X)
&\leqslant\sum_{n\geqslant1}\varPsi\Big(\frac{n}{X}\Big)\mu^2(n)|\lambda_f(n)|(\alpha+\beta)^{\omega(n)}\Big(\sum_{d|(n,P(z))}\varrho_d\Big)^2.\end{align*}
Taking $\alpha,\beta>0$ such that 
\begin{align}\label{eq:alpha-beta2}
\alpha+\beta\leqslant2,
\end{align}
so that
\begin{align*}H_{21}(X)
&\leqslant\sum_{n\geqslant1}\varPsi\Big(\frac{n}{X}\Big)\mu^2(n)|\lambda_f(n)|2^{\omega(n)}\Big(\sum_{d|(n,P(z))}\varrho_d\Big)^2.\end{align*}
By Cauchy's inequality, we have
\begin{align}\label{eq:H21(X)-Cauchy}
H_{21}(X)&\leqslant \sqrt{H_{21}'(X)H_{21}''(X)}\end{align}
with
\begin{align*}
H_{21}'(X)&=\sum_{n\geqslant1}\varPsi\Big(\frac{n}{X}\Big)\mu^2(n)|\lambda_f(n)|^22^{\omega(n)}\Big(\sum_{d|(n,P(z))}\varrho_d\Big)^2,\\
H_{21}''(X)&=\sum_{n\geqslant1}\varPsi\Big(\frac{n}{X}\Big)\mu^2(n)2^{\omega(n)}\Big(\sum_{d|(n,P(z))}\varrho_d\Big)^2.
\end{align*}

\subsection{Concluding an upper bound for $H_2(X)$}
The evaluations for $H_{21}'(X),H_{21}''(X)$ and $H_{22}(X)$ will rely on asymptotic computations of the average of Selberg sieve weights against some general multiplicative functions.
The later should be of independent interests and we will state a general situation by Theorem \ref{thm:asymptoticSelberg} in the appendix.

To evaluate $H_{22}(X)$ and $H_{21}''(X),$ we may take $h(n)=2^{\omega(n)}$ in Theorem \ref{thm:asymptoticSelberg}, so that
\begin{align}
H_{22}(X)&\leqslant(1+o(1))\fS\Big(\vartheta,\frac{\log X}{4\log z}\Big)\frac{X}{\log X}\label{eq:H22(X)-upperbound}\\
H_{21}''(X)&\leqslant(1+o(1))\fS\Big(\vartheta,\frac{\log X}{4\log z}\Big)\frac{X}{\log X},\label{eq:H21''(X)-upperbound}\end{align}
where $\fS(\cdot,\cdot)$ is given by \eqref{eq:fS(gamma,tau)}.

The evaluation of $H_{21}'(X)$ can be done by taking $h(n)=|\lambda_f(n)|^22^{\omega(n)}$ in Theorem \ref{thm:asymptoticSelberg}, and it suffices to verify 
the conditions of non-vanishing, meromorphic continuation \eqref{eq:meromorphiccontinuation},
first moment \eqref{eq:convergence-firstmoment} and second moment \eqref{eq:convergence-secondmoment} with some constants $L,c_0>0.$

In fact, it is well-known (see \cite[Proposition 2.3]{RS96} for instance) that
\[\sum_{p\leqslant x}\frac{\lambda_f(p)^2\log p}{p}=\log x+O_f(1),\]
which yields \eqref{eq:convergence-firstmoment} with some $L$ depending only on $f$. 
To check the condition \eqref{eq:meromorphiccontinuation} on meromorphic continuation, it may appeal to Lemma \ref{lm:pi-kappa(X)-moments} and derive that
\[\sum_{n\geqslant1}\mu^2(n)|\lambda_f(n)|^22^{\omega(n)}n^{-s}=\zeta(s)^2L(\sym^2f,s)^2F(s)\]
for $\Re s>1$ and $F(s)$ admits a Dirichlet series convergent absolutely in $\Re s>0.9.$ Hence the meromorphic continuation condition  \eqref{eq:meromorphiccontinuation} holds with $\cH^*(s)=L(\sym^2f,s)^2F(s)$ and $c_0=0.1$. The non-vanishing condition is guaranteed by the zero-free region of $L(\sym^2f,s)$ (see \cite[Theorem 5.44]{IK04} for instance).
After checking all above conditions, we conclude from Theorem \ref{thm:asymptoticSelberg} that 
\begin{align}
H_{21}'(X)&\leqslant(1+o(1))\fS\Big(\vartheta,\frac{\log X}{4\log z}\Big)\frac{X}{\log X}.\label{eq:H21'(X)-upperbound}\end{align}

In conclusion, Proposition \ref{prop:H2(X)-upperbound} follows immediately by combining \eqref{eq:H2(X)-upperbound-initial}, \eqref{eq:H21(X)-Cauchy}, \eqref{eq:H22(X)-upperbound}, \eqref{eq:H21''(X)-upperbound}, \eqref{eq:H21'(X)-upperbound}.

\smallskip

\section{Estimate for $H_3(X)$}
We rewrite $H_3(X)$ by
\begin{align*}
H_3(X)&=\sum_{n\geqslant1}\varPsi\Big(\frac{n}{X}\Big)\mu^2(n)(\lambda_f(n)-\eta\cdot\kl(1,n))\Big(\sum_{d|(n,P(z))}\varrho_d\Big)^2\\
&=\sum_{\substack{d\leqslant D\\d\mid P(z)}}\mu^2(d)\xi(d)\sum_{n\equiv0\bmod d}\varPsi\Big(\frac{n}{X}\Big)\mu^2(n)(\lambda_f(n)-\eta\cdot\kl(1,n))\end{align*}
with $\xi$ given by \eqref{eq:xi}.
In view of $|\xi(d)|\leqslant 3^{\omega(d)}$ for all squarefree $d\geqslant1$, Proposition \ref{prop:H3(X)-estimate} then follows from the following two lemmas.

\begin{lemma}\label{lm:BV-FM}
For any $A>0$, there exists some $B=B(A)>0$ such that
\begin{align*}
\sum_{q\leqslant \sqrt{X}\mathcal{L}^{-B}}3^{\omega(q)}\left|\sum_{n\equiv0\bmod q}\varPsi\Big(\frac{n}{X}\Big)\mu^2(n)\kl(1,n)\right|\ll X\mathcal{L}^{-A},\end{align*}
where the implied constant depends on $A$ and $\varPsi.$
\end{lemma}

\begin{lemma}\label{lm:EH-eigenvalues}
For any $A>0$, there exists some $B=B(A)>0$ such that
\begin{align*}
\sum_{q\leqslant X\mathcal{L}^{-B}}3^{\omega(q)}\left|\sum_{n\equiv0\bmod q}\varPsi\Big(\frac{n}{X}\Big)\mu^2(n)\lambda_f(n)\right|\ll X\mathcal{L}^{-A},\end{align*}
where the implied constant depends on $A,f$ and $\varPsi.$
\end{lemma}

Lemma \ref{lm:BV-FM}, which can be regarded as a Bombieri--Vinogradov type equidistribution for Kloosterman
sums, was initiated by Fouvry and Michel \cite{FM07} deriving from the spectral theory of automorphic forms without the weights $3^{\omega(q)}$ and $\mu^2(n)$. The current version is given by Sivak-Fischler \cite{SF09} and the author \cite{Xi15} with minor efforts.

Lemma \ref{lm:EH-eigenvalues} is not surprising to those readers that are familiar with automorphic forms, but the rigorous proof would require several extra lines. To simply the arguments, we assume the form $f$ is of level 1. In fact, the inner sum over $n$, denoted by $T$, can be rewritten as
\begin{align*}
T&=\sum_{d\leqslant2\sqrt{X}}\mu(d)
\sum_{n\equiv 0\bmod{[q,d^2]}}\varPsi\Big(\frac{n}{X}\Big)\lambda_f(n)\\
&=\sum_{d\leqslant2\sqrt{X}}\mu(d)
\sum_{n\geqslant1}\varPsi\Big(\frac{n[q,d^2]}{X}\Big)\lambda_f(n[q,d^2]).\end{align*}
By Hecke relation (see e.g., \cite[Formula (8.37)]{Iw02})
\[\lambda_f(mn)=\sum_{\ell\mid (m,n)}\mu(\ell)\lambda_f(m/\ell)\lambda_f(n/\ell),\]
we get
\begin{align*}
T
&=\sum_{d\leqslant2\sqrt{X}}\mu(d)\sum_{\ell\mid [q,d^2]}\mu(\ell)\lambda_f([q,d^2]/\ell)
\sum_{n\geqslant1}\varPsi\Big(\frac{n\ell [q,d^2]}{X}\Big)\lambda_f(n).\end{align*}
By partial summation and the well-known estimate (see e.g., \cite[Theorem 8.1]{Iw02})
\begin{align*}
\sum_{n\leqslant N}\lambda_f(n)\ll_f N^{\frac{1}{2}}\log N,\end{align*}
we derive that
\begin{align*}
T
&\ll_{f,g}\sqrt{X}\cL\sum_{d\leqslant2\sqrt{X}}\sum_{\ell\mid [q,d^2]}\mu^2(\ell)\frac{|\lambda_f([q,d^2]/\ell)|}{\sqrt{[q,d^2]\ell}}
\leqslant \sqrt{X}\cL\sum_{d\leqslant2\sqrt{X}}\frac{1}{[q,d^2]}\sum_{\ell\mid [q,d^2]}\mu^2(\ell)|\lambda_f(\ell)|\sqrt{\ell}.\end{align*}
Hence the original double sum in the lemma is bounded by
\begin{align*}
&\ll\sqrt{X}\cL\sum_{q\leqslant X\mathcal{L}^{-B}}3^{\omega(q)}\sum_{d\leqslant2\sqrt{X}}\frac{1}{[q,d^2]}\sum_{\ell\mid [q,d^2]}\mu^2(\ell)|\lambda_f(\ell)|\sqrt{\ell}\\
&\ll\sqrt{X}\cL\sum_{\ell\leqslant 4X\cL^{-B}}\mu^2(\ell)|\lambda_f(\ell)|\sqrt{\ell}\sum_{d\leqslant2\sqrt{X}}\frac{1}{d^2}
\sum_{\substack{q\leqslant X\mathcal{L}^{-B}\\ q\equiv0\bmod{\ell/(\ell,d^2)}}}\frac{3^{\omega(q)}(q,d^2)}{q}\\
&\ll \sqrt{X}\cL\sum_{\ell\leqslant4X\cL^{-B}}\frac{\mu^2(\ell)|\lambda_f(\ell)|}{\sqrt{\ell}}\sum_{d\leqslant2\sqrt{X}}\frac{3^{\omega(d)}(\ell,d)}{d^2}.\end{align*}
The lemma then follows from Cauchy's inequality and the Rankin--Selberg bound
\[\sum_{\ell\leqslant L}|\lambda_f(\ell)|^2\ll L,\]
as well as the choice $B=2A+4.$

\smallskip

\section{Numerical computations: concluding Theorems \ref{thm:non-identity} and \ref{thm:non-identity-general}}\label{sec:numerical}
In view of Propositions \ref{prop:H1(X)-lowerbound}, \ref{prop:H2(X)-upperbound} and \ref{prop:H3(X)-estimate},  we may conclude that
\begin{align}\label{eq:H(X)-positivelowerbound}
H^\pm(X)>\varepsilon_0 X\cL^{-1},
\end{align}
with some absolute constant $\varepsilon_0>0$, from the inequality
\begin{align}\label{eq:inequality-determingrho}
\rho\cdot \fA_1(\eta)>\fA_2(\eta)\end{align}
by choosing $\rho,\vartheta,z$ appropriately for a given $\eta\in\bR$, where
\begin{align}
\fA_1(\eta)&:=\sum_{2\leqslant i\leqslant7}I_i\cdot\sqrt{l_i^3/{u_i}},\nonumber\\
\fA_2(\eta)&:=(2+|\eta|)\fS\Big(\frac{1}{4},6\Big)=(2+|\eta|)16\ue^{2\gamma}\Big(\frac{2c_1(6)}{3}+\frac{c_2(6)}{9}\Big)\label{eq:A2(eta)}
\end{align}
subject to the restrictions \eqref{eq:alpha-beta1}, \eqref{eq:alpha-beta2} and the choice 
\begin{align}\label{eq:paremeter-choices}
\vartheta=\frac{1}{4},\ \ z=X^{\frac{1}{12}}.
\end{align}

\subsection{Upper bound for $\fA_2(\eta)$}
From the definitions \eqref{eq:sigma-equationsolution} and \eqref{eq:ff-equationsolution}, we find
\begin{align*}
\sigma(s)=
\begin{cases}
\dfrac{s^2}{8\ue^{2\gamma}},\ \ &s\in~]0,2],\\\noalign{\vskip 0,9mm}
\dfrac{s^2}{8\ue^{2\gamma}}\Big(4+\log4-2\log s-\dfrac{8s-4}{s^2}\Big),&s\in~]2,4],\\\noalign{\vskip 0,9mm}
\dfrac{s^2}{8\ue^{2\gamma}}\Big(4\displaystyle\int_4^s\dfrac{(t-2)^2\log(t-2)}{t^3}\ud t-(8+2\log4)\log s\\\noalign{\vskip 0,9mm}
\ \ \ \ \ \  +\dfrac{49+35\log4+8(\log4)^2}{4}-\dfrac{48+8\log4}{s}+\dfrac{32+4\log4}{s^2}\Big),&s\in~]4,6],
\end{cases}
\end{align*}
and
\begin{align*}
\ff(s)=
\begin{cases}
2,\ \ &s\in~]0,2],\\\noalign{\vskip 0,9mm}
4\log(s/2)+2,&s\in~]2,4],\\\noalign{\vskip 0,9mm}
8\displaystyle\int_4^s\dfrac{\log(t-2)}{t}\ud t-(8\log2-4)\log s+16(\log2)^2-4\log2+2,&s\in~]4,6].
\end{cases}
\end{align*}

Note that
\begin{align*}
c_1(6)&=\frac{1}{6}\int_0^6\sigma'(6-u)\ff(u)^2\ud u.
\end{align*}
From the positivity of $\sigma'$ and the monotonicity of $\ff$, it follows that
\begin{align*}
c_1(6)&=\frac{1}{6}\sum_{1\leqslant j\leqslant 6}\int_{j-1}^j\sigma'(6-u)\ff(u)^2\ud u\leqslant \frac{1}{6}\sum_{1\leqslant j\leqslant 6}\ff(j)^2\int_{j-1}^j\sigma'(6-u)\ud u\\
&=\frac{1}{6}\sum_{1\leqslant j\leqslant 6}\ff(j)^2(\sigma(7-j)-\sigma(6-j))\\
&=\frac{1}{6}\ff(1)^2\sigma(6)+\frac{1}{6}\sum_{3\leqslant j\leqslant 6}(\ff(j)^2-\ff(j-1)^2)\sigma(7-j).
\end{align*}

On the other hand,
\begin{align*}
c_2(6)&=\int_0^1\sigma'(6(1-u))\ud u\int_0^{3u}\ff(6u-2v)\{2\ff(6u)-\ff(6u-2v)\}\ud v\\
&=\frac{1}{12}\int_0^6\sigma'(6-u)\ud u\int_0^u\ff(v)\{2\ff(u)-\ff(v)\}\ud v.\end{align*}
Note that $\ff(v)\{2\ff(u)-\ff(v)\}\leqslant\ff(u)^2$ for all $v\in[0,u]$. Hence
\begin{align*}
c_2(6)&\leqslant\frac{1}{12}\int_0^6\sigma'(6-u)\ff(u)^2u\ud u.\end{align*}
From the positivity of $\sigma'$ and the monotonicity of $\ff$, it follows that
\begin{align*}
c_2(6)
&\leqslant\frac{1}{12}\sum_{1\leqslant j\leqslant6}\ff(j)^2j\int_{j-1}^{j}\sigma'(6-u)\ud u\\
&=\frac{1}{12}\sum_{1\leqslant j\leqslant6}\ff(j)^2j(\sigma(7-j)-\sigma(6-j))\\
&=\frac{1}{12}\ff(1)^2\sigma(6)+\frac{1}{12}\sum_{2\leqslant j\leqslant6}\{\ff(j)^2j-\ff(j-1)^2(j-1)\}\sigma(7-j).\end{align*}

Inserting the special values for $\sigma$ and $\ff$, we obtain
\begin{align*}
c_1(6)&\leqslant2.43762,\ \ \ \ c_2(6)\leqslant5.15051\end{align*}
upon the choice \eqref{eq:paremeter-choices}.
Combining the above two bounds and \eqref{eq:A2(eta)}, we conclude that
\begin{align*}
\fA_2(\eta)\leqslant 111.53(2+|\eta|).\end{align*}

\subsection{Lower bound for $\fA_1(\eta)$ and concluding Theorem \ref{thm:non-identity}}
With the help of Mathematica 10, we can find
\begin{align*}
I_2\geqslant 0.28768,\ \ \ \ & I_5\geqslant 0.14893\\
I_3\geqslant 1.04781,\ \ \ \ & I_6\geqslant 0.00424\\
I_4\geqslant 0.85019,\ \ \ \ & I_7\geqslant 7.25032\times10^{-6}.\end{align*}

For $\eta=\pm1$, we obtain
$\fA_1(\eta)\approx 3.687\times10^{-11},$ $\fA_2(\eta)\leqslant 334.59$, so that 
\eqref{eq:inequality-determingrho} holds by taking $\rho=9.076\times10^{12}$.
It suffices to solve the inequality 
\begin{align}\label{eq:tau-upperbound}
\tau_{\varDelta}(n;\alpha,\beta)<9.076\times10^{12}.\end{align}
To conclude Theorem \ref{thm:non-identity}, we should explore a lower bound for 
$\tau_\varDelta(n;\alpha,\beta)$, which grows as long as $\omega(n)$ increases.

Recall the definition \eqref{eq:tau(n;alpha,beta)} of the truncated divisor function $\tau_{\varDelta}(n;\alpha,\beta)$:
\begin{align*}
\tau_{\varDelta}(n;\alpha,\beta)=\sum_{\substack{d\mid n\\d\leqslant n^{\frac{1}{1+\varDelta}}}}\alpha^{\omega(d)}\beta^{\omega(n/d)}.\end{align*}
We would like to prove a lower bound for $\tau_{\varDelta}(n;\alpha,\beta)$ by elementary methods. To this end, let us recall a previous result of Soundararajan \cite{So92},
which gives a lower bound for the truncated convolution of multiplicative functions by {\it complete} convolutions.
The following lemma can be found in \cite[Theorem 4]{So92} with minor modifications on notation.

\begin{lemma}\label{lm:Soundararajan}
Let $t>0$ be a rational number and $g$ a multiplicative function with $0<g(p)\leqslant 1/t$ for all primes $p.$ Then, for each squarefree number $n\geqslant2,$ we have
\begin{align*}
\sum_{\substack{d\mid n\\ d\leqslant n^\frac{1}{1+t}}}g(d)\geqslant\fA(t)\sum_{d\mid n}g(d),
\end{align*}
where, if $t$ has the continued fraction expansion $[a_0,a_1,\cdots,a_k],$
\begin{align}\label{eq:continuedfraction}
\fA(t):=\frac{1}{1+a_0+a_1+\cdots+a_k}.
\end{align}
In particular, if $t$ is a positive integer, then $\fA(t)=1/(1+t).$
\end{lemma}

We now produce a lower bound for $\tau_{\varDelta}(n;\alpha,\beta)$ by virtue of Lemma \ref{lm:Soundararajan} subject to the restrictions \eqref{eq:alpha-beta1} and \eqref{eq:alpha-beta2}. Taking $\alpha,\beta,\varDelta$ such that $\alpha\varDelta=\beta>0,\varDelta\in\bQ~\cap~]1,+\infty[$, we conclude from Lemma \ref{lm:Soundararajan} that
\begin{align*}
\tau_{\varDelta}(n;\alpha,\beta)=\beta^{\omega(n)}\sum_{\substack{d\mid n\\d\leqslant n^{\frac{1}{1+\varDelta}}}}\Big(\frac{1}{\varDelta}\Big)^{\omega(d)}\geqslant \beta^{\omega(n)}\fA(\varDelta)\sum_{d\mid n}\Big(\frac{1}{\varDelta}\Big)^{\omega(d)}=\fA(\varDelta)(\alpha+\beta)^{\omega(n)}.\end{align*}

Following the above arguments, we are now in a position to solve the inequality
\begin{align*}
\fA(\varDelta)(\alpha+\beta)^{\omega(n)}<9.076\times10^{12},\end{align*}
where $\varDelta,\alpha,\beta>0$ are chosen freely subject to the following restrictions
\begin{align*}
\varDelta=\beta/\alpha\in\bQ~\cap~]1,+\infty[,\ \ \ \ \frac{8\alpha}{3\pi}+2\beta\leqslant2,\ \ \ \ \alpha+\beta\leqslant2.
\end{align*}
In particular, we would like to take 
\begin{align*}
\varDelta=\frac{14}{13},\ \ \ \ \alpha=\frac{39\pi}{52+42\pi},\ \ \ \ \beta=\frac{21\pi}{26+21\pi},
\end{align*}
in which case one has $\fA(\varDelta)=\frac{1}{15}.$ It now suffices to 
solve the inequality
\begin{align*}
\frac{1}{15}\Big(\frac{81\pi}{52+42\pi}\Big)^{\omega(n)}<9.076\times10^{12},\end{align*}
which yields $\omega(n)<100.29,$ i.e., $\omega(n)\leqslant100.$

To conclude the quantitative statement in Theorem \ref{thm:non-identity}, we would like to argue as follows.
Put $\cN(X):=\{n\in[X,2X]:\lambda_f(n)>\kl(1,n),\omega(n)\leqslant 100,\mu^2(n)=1\}$. Trivially, we have
\begin{align*}
H^+(X)
&\leqslant\rho\sum_{\tau_{\varDelta}(n;\alpha,\beta)<\rho}\varPsi\Big(\frac{n}{X}\Big)\mu^2(n)\{|\psi(n)|+\psi(n)\}\Big(\sum_{d|(n,P(z))}\varrho_d\Big)^2\\
&\leqslant 2\rho\sum_{\substack{\psi(n)>0\\ \omega(n)\leqslant 100}}\varPsi\Big(\frac{n}{X}\Big)\mu^2(n)|\psi(n)|\Big(\sum_{d|(n,P(z))}\varrho_d\Big)^2\end{align*}
with $\psi(n)=\lambda_f(n)-\kl(1,n).$
By Cauchy's inequality, we find
\begin{align*}
H^+(X)^2
&\leqslant4\rho^2|\cN(X)|\sum_{\substack{\psi(n)>0\\ \omega(n)\leqslant 100}}\varPsi^2\Big(\frac{n}{X}\Big)\mu^2(n)|\psi(n)|^2\Big(\sum_{d|(n,P(z))}\varrho_d\Big)^4.\end{align*}
Note that
\begin{align*}
\Bigg|\sum_{d|(n,P(z))}\varrho_d\Bigg|\leqslant 2^{\omega(n)}\end{align*}
for each squarefree $n$, from which and Weil's bound for Kloosterman sums, it follows that
\begin{align*}
H^+(X)^2
&\leqslant4\rho^2|\cN(X)|\sum_{\omega(n)\leqslant 100}\varPsi^2\Big(\frac{n}{X}\Big)\mu^2(n)|\psi(n)|^24^{\omega(n)}\Big(\sum_{d|(n,P(z))}\varrho_d\Big)^2\\
&\leqslant4^{101}\rho^2|\cN(X)|\sum_{\omega(n)\leqslant 100}\varPsi^2\Big(\frac{n}{X}\Big)\mu^2(n)(|\lambda_f(n)|^2+4^{100})\Big(\sum_{d|(n,P(z))}\varrho_d\Big)^2\\
&\leqslant4^{101}\rho^2|\cN(X)|\sum_{n\geqslant1}\varPsi^2\Big(\frac{n}{X}\Big)\mu^2(n)(|\lambda_f(n)|^2+4^{100})\Big(\sum_{d|(n,P(z))}\varrho_d\Big)^2.\end{align*}
We now proceed as in the proof of Proposition \ref{prop:H2(X)-upperbound}, and the last sum over $n$ can be bounded by $O(X\cL^{-1})$ with an absolute constant. Therefore,
\begin{align*}
H^+(X)^2\ll X\cL^{-1}\cdot|\cN(X)|.\end{align*}
Combining this with \eqref{eq:H(X)-positivelowerbound}, we then arrive at
\begin{align*}
|\cN(X)|\gg X\cL^{-1}.\end{align*}

Similar arguments can also lead to 
\begin{align*}
|\{n\in[X,2X]:\lambda_f(n)<\kl(1,n),\omega(n)\leqslant 100,\mu^2(n)=1\}|\gg X\cL^{-1}.\end{align*}
We now complete the proof of Theorem \ref{thm:non-identity}.

\subsection{The case of general $\eta$}
Given an $\eta\in\bR$, one may see that those $l_i$'s in Proposition \ref{prop:Sigma(X,alphai)-lowerbound} are not always positive. To obtain a positive lower bound for $\fA_1(\eta)$, we need to solve the inequality $\fA_1(\eta)>0$, which holds provided that
\begin{align}\label{eq:A1(eta)positive-eta}
|\eta|\in[0,~1.23]\cup[11.84,~+\infty[.
\end{align}
For such $\eta$ we may choose a considerably large $\rho$ such that \eqref{eq:inequality-determingrho} holds, and thus we can always produce almost primes in Theorem \ref{thm:non-identity-general} for a general $\eta\in\bR$ satisfying \eqref{eq:A1(eta)positive-eta}. 

In fact, as $|\eta|$ is sufficiently large, we find from Lemma \ref{prop:Sigma(X,alphai)-lowerbound}  that
$$\fA_1(\eta)\geqslant c_1|\eta|,\ \ \ \fA_2(\eta)\leqslant c_2|\eta|$$
for some constant $c_1,c_2>0.$ Therefore, a certain absolute $\rho$ could be found for all such large $|\eta|$, for which we may explore a uniform $r$ in 
Theorem \ref{thm:non-identity-general}.
This is not surprising since Kloosterman sums will dominate the contributions to $H^\pm(X)$ if $|\eta|$ is quite large, and the difficulty of Theorem 
\ref{thm:non-identity-general} becomes close to the sign changes of Kloosterman sums with almost prime moduli, as considered in \cite{FM03b,FM07,SF09,Ma11,Xi15,Xi18}.

On the other hand, if $|\eta|$ decays to zero, we also have uniform bounds for $\fA_1(\eta)$ and $\fA_2(\eta)$. Following a similar argument,
the choice of $r$ in Theorem \ref{thm:non-identity-general} can also be made uniformly
in all such small $|\eta|$.

It remains to consider the complementary range of $\eta$ to \eqref{eq:A1(eta)positive-eta}. Recall that
\begin{align*}l_i=(1-4|\eta|\cdot(\tfrac{8}{3\pi})^{i-1}(\tfrac{11}{12})^i+B_i\eta^2)^+\end{align*}
in Proposition \ref{prop:Sigma(X,alphai)-lowerbound},
and the positivity of $l_i$ lies in the essential part of this paper. For any $\eta$ with
$|\eta|\in[1.23,~11.94],$ one may find $l_i>0.2$ as long as $i\geqslant 17.$ Therefore, one may sum up to $i=17$ in \eqref{eq:H1(X)-lowerbound-initial} with
\begin{align*}
\mathcal{R}_i&:=\{(\alpha_2,\cdots,\alpha_i)\in[\tfrac{1}{18},1[^{i-1}:\tfrac{1}{2}(1-\alpha_2-\cdots-\alpha_i)<\alpha_2+\cdots+\alpha_{i-1}\}\\
&\ \ \ \ \ \cap\{(\alpha_2,\cdots,\alpha_i)\in[\tfrac{1}{18},1[^{i-1}:\tfrac{1}{2}(\alpha_3+\cdots+\alpha_i)<\alpha_2\}\\
&\ \ \ \ \ \cap\{(\alpha_2,\cdots,\alpha_i)\in[\tfrac{1}{18},1[^{i-1}:\alpha_i<\alpha_{i-1}<\cdots<\alpha_2<1-\alpha_2-\cdots-\alpha_i\}
\end{align*}
for $i=17$. To evaluate the Selberg sieve weight, we may re-take $z=X^{\frac{1}{19}}$, so that
\begin{align*}\sum_{d|(n,P(z))}\varrho_d=\varrho_1=1\end{align*}
if $n$ is restricted to $\cP_{17}(X,\boldsymbol\alpha_{17})$.
Following the above arguments in proving Proposition \ref{prop:H1(X)-lowerbound}, we may obtain a positive lower bound for $H_{1,17}(X)$, and thus that for $H_1(X).$
To complete the proof of Theorem \ref{thm:non-identity-general}, it remains to produce an explicit numerical upper bound for $H_2(X).$
This requires a delicate analysis on $\sigma(s),\ff(s)$, and the details are omitted here.

The Mathematica codes can be found at \url{http://gr.xjtu.edu.cn/web/ping.xi/miscellanea} or requested from the author.

\smallskip

\appendix

\section{Multiplicative functions against M\"obius}

We would like to evaluate a weighted average of general multiplicative functions against M\"obius function.
This will be employed in the evaluation of Selberg sieve weights essentially given by \eqref{eq:varrho_d}.

Let $g$ be a non-negative multiplicative function with $0\leqslant g(p)<1$ for each $p\in\cP$. Suppose the Dirichlet series
\begin{align}\label{eq:G(s)}
\cG(s):=\sum_{n\geqslant1}\mu^2(n)g(n)n^{-s}
\end{align}
converges absolutely for $\Re s>1.$ Assume there exist a positive integer $\kappa$ and some constants $L,c_0>0,$ such that
\begin{align}\label{eq:meromorphiccontinuation-G}
\cG(s)=\zeta(s+1)^\kappa\cF(s),
\end{align}
where $\cF(s)$ is holomorphic for $\Re s\geqslant -c_0$ and does not vanish in the region
\begin{align}\label{eq:zerofree}
\cD:=\Big\{\sigma+it: t\in\bR,\sigma\geqslant-\frac{1}{L\cdot\log(|t|+2)}\Big\},
\end{align}
and $|1/\cF(s)|\leqslant L$ for all $s\in\cD.$
We also assume
\begin{align}\label{eq:convergence-firstmoment-g}
\left|\sum_{p\leqslant x}g(p)\log p-\kappa\log x\right|\leqslant L
\end{align}
holds for all $x\geqslant3$ and
\begin{align}\label{eq:convergence-secondmoment-g}
\sum_{p}g(p)^2p^{2c_0}<+\infty.
\end{align}

We are interested in the asymptotic behaviour of the sum
\begin{align*}
\cM_\kappa(x,z;q)&=\sum_{\substack{n\leqslant x\\ n\mid P(z)\\ (n,q)=1}}\mu(n)g(n)\Big(\log\frac{x}{n}\Big)^\kappa,\end{align*}
where $q$ is a positive integer and $x,z\geqslant3$.

\begin{lemma}\label{lm:averageagainstMobius}
Let $q\geqslant1.$ Under the assumption as above, we have
\begin{align*}
\cM_\kappa(x,z;q)&=H\cdot\prod_{p\mid q}(1-g(p))^{-1}\cdot m_\kappa(s)+O(\kappa^{\omega(q)}(\log z)^{-A})\end{align*}
for all $A>0,x\geqslant2,z\geqslant2$ with $x\leqslant z^{O(1)}$, where $s=\log x/\log z,$
\[H=\prod_{p}(1-g(p))\Big(1-\frac{1}{p}\Big)^{-\kappa},\]
and $m_\kappa(s)$ is a continuous solution to the differential-difference equation
\begin{align}\label{eq:mkappa(s)}
\begin{cases}
m_\kappa(s)=\kappa!,\ \ &s\in~]0,1],\\
sm_\kappa'(s)=\kappa m_\kappa(s-1),&s\in~]1,+\infty[.\end{cases}
\end{align}
The implied constant depends on $A,\kappa,L$ and $c_0.$
\end{lemma}

\proof
We are inspired by \cite[Appendix A.3]{FI10}.
Write $\cM_\kappa(x,x;q)=\cM_\kappa(x;q)$. By Mellin inversion, we have
\begin{align*}
\cM_\kappa(x;q)&=\sum_{\substack{n\leqslant x\\ (n,q)=1}}\mu(n)g(n)\Big(\log\frac{x}{n}\Big)^\kappa=\frac{\kappa!}{2\pi i}\int_{2-i\infty}^{2+i\infty}\cG(t,q)\frac{x^t}{t^{\kappa+1}}\ud t,\end{align*}
where
\begin{align*}
\cG(t,q)=\sum_{\substack{n\geqslant1\\(n,q)=1}}\frac{\mu(n)g(n)}{n^t},\ \ \Re t>1.\end{align*}
Note that
\begin{align*}
\cG(t,q)=\prod_{p\nmid q}\Big(1-\frac{g(p)}{p^t}\Big)=\prod_{p\mid q}\Big(1-\frac{g(p)}{p^t}\Big)^{-1}\frac{\cG^*(t)}{\zeta(t+1)^\kappa},\end{align*}
where
\begin{align*}
\cG^*(t)=\prod_{p}\Big(1-\frac{g(p)}{p^t}\Big)\Big(1-\frac{1}{p^{t+1}}\Big)^{-\kappa}=\prod_{p}\Big(1-\frac{g(p)^2}{p^{2t}}\Big)\frac{1}{\cF(t)},\end{align*}
which is absolutely convergent and holomorphic for $t\in\cC$ by \eqref{eq:meromorphiccontinuation-G}, \eqref{eq:convergence-firstmoment-g} and \eqref{eq:convergence-secondmoment-g}.
Hence we find
\begin{align*}
\cM_\kappa(x;q)
&=\frac{\kappa!}{2\pi i}\int_{2-i\infty}^{2+i\infty}\prod_{p\mid q}\Big(1-\frac{g(p)}{p^t}\Big)^{-1}\frac{\cG^*(t)x^t}{\zeta(t+1)^\kappa t^{\kappa+1}}\ud t.\end{align*}
Shifting the $t$-contour to the left boundary of $\cC$ and passing one simple pole at $t=0$, we get
\begin{align*}
\cM_\kappa(x;q)&=\kappa!\cG^*(0)\prod_{p\mid q}(1-g(p))^{-1}+O(\kappa^{\omega(q)}(\log 2x)^{-A})\end{align*}
for any fixed $A>0$.

For $s=\log x/\log z,$  we expect that
\begin{align}\label{eq:expectation-generalz}
\cM_\kappa(x,z;q)&=c(q)m_\kappa(s)+O(\kappa^{\omega(q)}(\log z)^{-A})\end{align}
for all $A>0,x\geqslant2,z\geqslant2$ and $q\geqslant1$, where $c(q)$ is some constant defined in terms of $g$ and depending also on $q$, and $m_\kappa(s)$ is a suitable continuous function in $s>0.$
As mentioned above, this expected asymptotic formula holds for $0<s\leqslant1,$ in which case we may take
\begin{align*}
c(q)=\cG^*(0)\prod_{p\mid q}(1-g(p))^{-1},\ \ \ m_\kappa(s)=\kappa!.\end{align*}
We now move to the case $s>1$ and prove the asymptotic formula \eqref{eq:expectation-generalz} by induction. Since $x\leqslant z^{O(1)},$ this induction will have a bounded number of steps.
We first consider the difference $\cM_\kappa(x,z;q)-\cM_\kappa(x;q)$. In fact, each $n$ that contributes to this difference has a prime factor at least $z$, and we may decompose $n=mp$ uniquely up to the restriction $z\leqslant p<x,$ $m\mid P(p).$ Hence
\begin{align}
\cM_\kappa(x,z;q)
&=\cM_\kappa(x;q)+\sum_{\substack{z\leqslant p<x\\(p,q)=1}}g(p)\sum_{\substack{m\leqslant x/p\\ m\mid P(p)\\ (m,q)=1}}\mu(m)g(m)\Big(\log\frac{x}{mp}\Big)^\kappa\nonumber\\
&=\cM_\kappa(x;q)+\sum_{\substack{z\leqslant p<x\\(p,q)=1}}g(p)\cM_\kappa(x/p,p;q).\label{eq:iteration}\end{align}

Substituting \eqref{eq:expectation-generalz} to \eqref{eq:iteration}, we get
\begin{align*}
\cM_\kappa(x,z;q)
&=c(q)\kappa!+c(q)\sum_{\substack{z\leqslant p<x\\(p,q)=1}}g(p)m_\kappa\Big(\frac{\log(x/p)}{\log p}\Big)+O(\kappa^{\omega(q)}(\log x)^{-A})\\
&\ \ \ \ +O\Big(\kappa^{\omega(q)}\sum_{\substack{z\leqslant p<x\\(p,q)=1}}g(p)(\log (2x/p))^{-A}\Big).\end{align*}
By partial summation, we find
\begin{align*}
\cM_\kappa(x,z;q)
&=c(q)\Big\{\kappa!+\kappa\int_1^sm_\kappa\Big(\frac{s}{u}-1\Big)\frac{\ud u}{u}\Big\}+O(\kappa^{\omega(q)}(\log z)^{-A}).\end{align*}
Hence, by \eqref{eq:expectation-generalz}, $m_\kappa(s)$ should satisfy the equation
\begin{align*}
m_\kappa(s)=\kappa!+\kappa\int_1^sm_\kappa\Big(\frac{s}{u}-1\Big)\frac{\ud u}{u}=\kappa!+\kappa\int_1^sm_\kappa(u-1)\frac{\ud u}{u}\end{align*}
for $s>1$.
Taking the derivative with respect to $s$ gives \eqref{eq:mkappa(s)}.
\endproof

\begin{remark}
To extend $m_\kappa(s)$ to be defined on $\bR,$ we may put $m_\kappa(s)=0$ for $s\leqslant0$.
\end{remark}

\bigskip

\section{A two-dimensional Selberg sieve with asymptotics}
This section devotes to present a two-dimensional Selberg sieve that plays an essential role in proving Proposition \ref{prop:H2(X)-upperbound}.

Let $h$ be a non-negative multiplicative function. Suppose the Dirichlet series
\begin{align}\label{eq:H(s)}
\cH(s):=\sum_{n\geqslant1}\mu^2(n)h(n)n^{-s}
\end{align}
converges absolutely for $\Re s>1.$ Assume there exist some constants $L,c_0>0,$ such that
\begin{align}\label{eq:meromorphiccontinuation}
\cH(s)=\zeta(s)^2\cH^*(s),
\end{align}
where $\cH^*(s)$ is holomorphic for $\Re s\geqslant 1-c_0,$ and does not vanish in the region
$\cD$ as given by \eqref{eq:zerofree} and $|1/\cH^*(s)|\leqslant L$ for all $s\in\cD$.
We also assume
\begin{align}\label{eq:convergence-firstmoment}
\left|\sum_{p\leqslant x}\frac{h(p)\log p}{p}-2\log x\right|\leqslant L
\end{align}
holds for all $x\geqslant3$ and
\begin{align}\label{eq:convergence-secondmoment}
\sum_{p}h(p)^2p^{2c_0-2}<+\infty.
\end{align}

Define
\begin{align*}
S(X,z;h,\boldsymbol\varrho)=\sum_{n\geqslant1}\varPsi\Big(\frac{n}{X}\Big)\mu^2(n)h(n)\Big(\sum_{d|(n,P(z))}\varrho_d\Big)^2,
\end{align*}
where $\boldsymbol\varrho=(\varrho_d)$ is given as in \eqref{eq:varrho_d} and $\varPsi$ is a fixed non-negative smooth function supported in $[1,2]$ with normalization \eqref{eq:normalization-Psi}.

\begin{theorem}\label{thm:asymptoticSelberg}
Let  $X,D,z\geqslant3$ with $X\leqslant D^{O(1)}$ and $X\leqslant z^{O(1)}.$ Put $\tau=\log D/\log z$ and $\sqrt{D}=X^\vartheta\exp(-\sqrt{\cL}),\vartheta\in~]0,\frac{1}{2}[.$
Under the above assumptions, we have
\begin{align*}
S(X,z;h,\boldsymbol\varrho)&=(1+o(1))\fS(\vartheta,\tau)X\cL^{-1},\end{align*}
where $\fS(\vartheta,\tau)$ is defined by 
\begin{align}\label{eq:fS(gamma,tau)}
\fS(\vartheta,\tau)&=16\ue^{2\gamma}\Big(\frac{c_1(\tau)}{4\tau\vartheta^2}+\frac{c_2(\tau)}{\tau^2\vartheta}\Big),\end{align}
where
\begin{align*}
c_1(\tau)&=\int_0^1\sigma'((1-u)\tau)\ff(u\tau)^2\ud u,\\
c_2(\tau)&=\int_0^1\int_0^1\sigma'((1-u)\tau)\ff(u\tau-2v)\{2\ff(u\tau)-\ff(u\tau-2v)\}\ud u\ud v.\end{align*}
Here $\sigma(s)$ is the continuous solution to the differential-difference equation
\begin{align}\label{eq:sigma-equationsolution}
\begin{cases}
\sigma(s)=\dfrac{s^2}{8\ue^{2\gamma}},\ \ &s\in~]0,2],\\\noalign{\vskip 0,9mm}
(s^{-2}\sigma(s))'=-2s^{-3}\sigma(s-2),&s\in~]2,+\infty[,
\end{cases}
\end{align}
and $\ff(s)=m_2(s/2)$ as given by $\eqref{eq:mkappa(s)},$ i.e., $\ff(s)$ is the continuous solution to the differential-difference equation
\begin{align}\label{eq:ff-equationsolution}
\begin{cases}
\ff(s)=2,\ \ &s\in~]0,2],\\
s\ff'(s)=2\ff(s-2),&s\in~]2,+\infty[.\end{cases}
\end{align} 
\end{theorem}

\begin{remark}
Theorem \ref{thm:asymptoticSelberg} is a generalization of \cite[Proposition 4.1]{Xi18} with a general multiplicative function $h$ and the extra restriction $d\mid P(z)$,
but specializing $k=2$ therein. It would be rather interesting to extend the case to a general $k\in\bZ^+$ and we would like to concentrate this problem in the near future.

We now choose $z=\sqrt{D}$, so that the restriction $d\mid P(z)$ is redundant, in which case one has $\tau=2.$ Note that
\begin{align*}
c_1(2)&=4\int_0^1\sigma'(2u)\ud u=\frac{1}{\ue^{2\gamma}},\\
c_2(2)&=4\int_0^1\sigma'(2(1-u))u\ud u=\frac{1}{3\ue^{2\gamma}}.\end{align*}
For $\vartheta=1/4,$ we find $\fS(\vartheta,\tau)=\fS(1/4,2)=112/3$, which coincides with $4\fc(2,F)$ in \cite[Proposition 4.1]{Xi18} by taking $F(x)=x^2$ therein.
\end{remark}

We now give the proof of Theorem \ref{thm:asymptoticSelberg}.
To begin with, we write by \eqref{eq:xi} that
\begin{align*}
S(X,z;h,\boldsymbol\varrho)
&=\sum_{d\mid P(z)}\xi(d)\sum_{n\equiv0\bmod d}\varPsi\Big(\frac{n}{X}\Big)\mu^2(n)h(n)\\
&=\sum_{d\mid P(z)}\xi(d)h(d)\sum_{(n,d)=1}\varPsi\Big(\frac{nd}{X}\Big)\mu^2(n)h(n).\end{align*}
By Mellin inversion,
\begin{align*}
\sum_{(n,d)=1}\varPsi\Big(\frac{nd}{X}\Big)\mu^2(n)h(n)
&=\frac{1}{2\pi i}\int_{(2)}\widetilde{\varPsi}(s)(X/d)^s\cH^\flat(s,d)\ud s,\end{align*}
where, for $\Re s>1,$
\[\cH^\flat(s,d)=\sum_{\substack{n\geqslant1\\(n,d)=1}}\frac{\mu^2(n)h(n)}{n^s}.\]

For $\Re s>1,$ we first write
\begin{align*}\cH^\flat(s,d)&=\prod_{p\nmid d}\Big(1+\frac{h(p)}{p^s}\Big)=\prod_{p\mid d}\Big(1+\frac{h(p)}{p^s}\Big)^{-1}\cH(s)=\prod_{p\mid d}\Big(1+\frac{h(p)}{p^s}\Big)^{-1}\zeta(s)^2\cG(s).\end{align*}
Note that
\begin{align*}
\cG(1)=\lim_{s\rightarrow1}\frac{\cH(s)}{\zeta(s)^2}=\prod_p\Big(1+\frac{h(p)}{p}\Big)\Big(1-\frac{1}{p}\Big)^2.\end{align*}
By \eqref{eq:meromorphiccontinuation}, $\cH^\flat(s,d)$ admits a meromorphic continuation to $\Re s\geqslant 1-c_0.$ Shifting the $s$-contour to the left beyond $\Re s=1,$ we may obtain
\begin{align*}
\sum_{(n,d)=1}\varPsi\Big(\frac{nd}{X}\Big)\mu^2(n)h(n)
&=\Res_{s=1}\widetilde{g}(s)\cG(s)(X/d)^s\prod_{p\mid d}\Big(1+\frac{h(p)}{p^s}\Big)^{-1}\zeta(s)^2+O((X/d)\cL^{-100}).\end{align*}
We compute the residue as
\begin{align*}
\Res_{s=1}[\cdots]
&=\frac{\ud}{\ud s}\widetilde{\varPsi}(s)\cG(s)(X/d)^s\prod_{p\mid d}\Big(1+\frac{h(p)}{p^s}\Big)^{-1}\zeta(s)^2(s-1)^2\Big|_{s=1}\\
&=\widetilde{\varPsi}(1)\cG(1)\prod_{p\mid d}\Big(1+\frac{h(p)}{p}\Big)^{-1}\frac{X}{d}\Big(\log(X/d)+\sum_{p\mid d}\frac{h(p)\log p}{p+h(p)}+c\Big)\\
&=\cG(1)\prod_{p\mid d}\Big(1+\frac{h(p)}{p}\Big)^{-1}\frac{X}{d}\Big(\log X-\sum_{p\mid d}\frac{p\log p}{p+h(p)}+c\Big),\end{align*}
where $c$ is some constant independent of $d$.

Define $\beta$ and $\beta^*$ to be multiplicative functions supported on squarefree numbers via
\begin{align*}\beta(p)=\frac{p}{h(p)}+1,\ \ \ \ \beta^*(p)=\beta(p)-1=\frac{p}{h(p)}.\end{align*}
Define $L$ to be an additive function supported on squarefree numbers via
\begin{align*}
L(p)=\frac{\beta^*(p)\log p}{\beta(p)}.\end{align*}
Therefore, for each squarefree number $d$, we have
\begin{align*}
\beta(d)=\prod_{p\mid d}\Big(\frac{p}{h(p)}+1\Big),\ \ \ \ \beta^*(d)=\frac{d}{h(d)},\ \ \ \ L(d)= \sum_{p\mid d}\frac{\beta^*(p)\log p}{\beta(p)}.\end{align*}
In this way, we may obtain
\begin{align*}S(X;h,\boldsymbol\varrho)
&=\cG(1)X\{S_1(X)\cdot (\log X+c)-S_2(X)\}+O(X\cL^{-2}),\end{align*}
where
\begin{align*}
S_1(X)
&=\sum_{d\mid P(z)}\frac{\xi(d)}{\beta(d)},\\
S_2(X)
&=\sum_{d\mid P(z)}\frac{\xi(d)}{\beta(d)}L(d).\end{align*}

Note that
\begin{align*}
S_1(X)&=\mathop{\sum\sum}_{d_1,d_2\mid P(z)}\frac{\varrho_{d_1}\varrho_{d_2}}{\beta([d_1,d_2])}\\
&=\mathop{\sum\sum}_{d_1,d_2\mid P(z)}\frac{\varrho_{d_1}\varrho_{d_2}}{\beta(d_1)\beta(d_2)}\beta((d_1,d_2))\\
&=\mathop{\sum\sum}_{d_1,d_2\mid P(z)}\frac{\varrho_{d_1}\varrho_{d_2}}{\beta(d_1)\beta(d_2)}\sum_{l\mid(d_1,d_2)}\beta^*(l).\end{align*}
Hence we may diagonalize $S_1(X)$ by
\begin{align}\label{eq:S1(X)-initial}
S_1(X)
&=\sum_{\substack{l\leqslant\sqrt{D}\\l\mid P(z)}}\beta^*(l)y_l^2,\end{align}
where, for each $l\mid P(z)$ and $l\leqslant\sqrt{D},$
\begin{align*}
y_l=\sum_{\substack{d\mid P(z)\\ d\equiv0\bmod l}}\frac{\varrho_d}{\beta(d)}.\end{align*}

From the definition of sieve weights \eqref{eq:varrho_d}, we find
\begin{align*}
y_l
&=\frac{4\mu(l)}{\beta(l)(\log D)^2}\sum_{\substack{d\leqslant\sqrt{D}/l\\dl\mid P(z)}}\frac{\mu(d)}{\beta(d)}\Big(\log\frac{\sqrt{D}/l}{d}\Big)^2.\end{align*}
Applying Lemma \ref{lm:averageagainstMobius} with $g(p)=1/\beta(p)$ and $q=l$, we have
\begin{align}\label{eq:yl}
y_l
&=\frac{4\mu(l)}{\cG(1)\beta^*(l)(\log D)^2}m_2\Big(\frac{\log(\sqrt{D}/l)}{\log z}\Big)+O\Big(\frac{\tau(l)}{\beta(l)}(\log z)^{-A}\Big).\end{align}
Inserting this expression to \eqref{eq:S1(X)-initial}, we have
\begin{align*}
S_1(X)
&=\frac{16(1+o(1))}{\cG(1)^2(\log D)^4}\sum_{\substack{l\leqslant\sqrt{D}\\l\mid P(z)}}\frac{1}{\beta^*(l)}m_2\Big(\frac{\log(\sqrt{D}/l)}{\log z}\Big)^2.\end{align*}

Following \cite[Lemma 6.1]{HR74}, we have
\begin{align}
\sum_{\substack{l\leqslant x\\l\mid P(z)}}\frac{1}{\beta^*(l)}=\frac{1}{W(z)}\Big\{\sigma(2\log x/\log z)+O\Big(\frac{(\log x/\log z)^5}{\log z}\Big)\Big\}\end{align}
with
\begin{align*}W(z)&=\prod_{p<z}\Big(1-\frac{1}{\beta(p)}\Big),\end{align*}
from which and partial summation, we find
\begin{align*}
S_1(X)
&=\frac{16\tau c_1(\tau)}{\cG(1)^2W(z)(\log D)^4}\cdot (1+o(1))\end{align*}
with $\tau=\log D/\log z$ and
\begin{align}
c_1(\tau)=\int_0^1\sigma'((1-u)\tau)\ff(u\tau)^2\ud u.\end{align}

We now turn to consider $S_2(X)$. Note that $L(d)$ is an additive function supported on squarefree numbers. We then have
\begin{align*}
S_2(X)
&=\mathop{\sum\sum}_{d_1,d_2\mid P(z)}\frac{\varrho_{d_1}\varrho_{d_2}}{\beta([d_1,d_2])}L([d_1,d_2])\\
&=\mathop{\sum\sum}_{dd_1d_2\mid P(z)}\frac{\varrho_{dd_1}\varrho_{dd_2}}{\beta(dd_1d_2)}\{L(d)+L(d_1)+L(d_2)\},\end{align*}
where there is an implicit restriction that $d,d_1,d_2$ are pairwise coprime. By M\"obius formula, we have
\begin{align*}
S_2(X)
&=\mathop{\sum\sum\sum}_{dd_1,dd_2\mid P(z)}\frac{\varrho_{dd_1}\varrho_{dd_2}}{\beta(d)\beta(d_1)\beta(d_2)}\{L(d)+L(d_1)+L(d_2)\}\sum_{l\mid(d_1,d_2)}\mu(l)\\
&=\mathop{\sum\sum\sum\sum}_{ldd_1,ldd_2\mid P(z)}\frac{\mu(l)\varrho_{ldd_1}\varrho_{ldd_2}}{\beta(l)^2\beta(d)\beta(d_1)\beta(d_2)}\{L(ldd_1)+L(ldd_2)-L(d)\}\\
&=2S_{21}(X)-S_{22}(X)\end{align*}
with
\begin{align*}
S_{21}(X)&=\sum_{l\mid P(z)}\beta^*(l)y_ly_l',\\
S_{22}(X)&=\sum_{l\mid P(z)}v(l)y_l^2,\end{align*}
where for each $l\mid P(z),l\leqslant\sqrt{D},$
\begin{align*}
y_l'=\sum_{\substack{d\mid P(z)\\ d\equiv0\bmod l}}\frac{\varrho_dL(d)}{\beta(d)}.\end{align*}
and
\begin{align}\label{eq:v(d)}
v(l)&=\beta(l)\sum_{uv=l}\frac{\mu(u)L(v)}{\beta(u)}.\end{align}
Moreover, we have
\begin{align*}
y_l'&=\sum_{dl\mid P(z)}\frac{\varrho_{dl}L(dl)}{\beta(dl)}=\sum_{d\mid P(z)}\frac{\varrho_{dl}L(d)}{\beta(dl)}+L(l)y_l\\
&=\sum_{p<z}\frac{\beta^*(p)\log p}{\beta(p)}\sum_{d\mid P(z)}\frac{\varrho_{pdl}}{\beta(pdl)}+L(l)y_l\\
&=\sum_{p<z}\frac{y_{pl}\beta^*(p)\log p}{\beta(p)}+L(l)y_l.\end{align*}
It then follows that
\begin{align*}
S_{21}(X)&=\sum_{p<z}\frac{\beta^*(p)\log p}{\beta(p)}\sum_{l\mid P(z)}\beta^*(l)y_ly_{pl}+\sum_{l\mid P(z)}L(l)\beta^*(l)y_l^2\\
&=\sum_{p<z}\frac{\beta^*(p)\log p}{\beta(p)}\sum_{pl\mid P(z)}\beta^*(l)y_ly_{pl}+\sum_{p<z}\frac{\beta^*(p)^2\log p}{\beta(p)}\sum_{pl\mid P(z)}\beta^*(l)y_{pl}^2\\
&=S_{21}'(X)+S_{21}''(X),\end{align*}
say.

From \eqref{eq:yl}, it follows, by partial summation, that
\begin{align*}
S_{21}'(X)
&=-\frac{16(1+o(1))}{\cG(1)^2(\log D)^4}\sum_{l\mid P(z)}\frac{1}{\beta^*(l)}m_2\Big(\frac{\log(\sqrt{D}/l)}{\log z}\Big)\sum_{\substack{p<z\\ p\nmid l}}\frac{\log p}{\beta(p)}m_2\Big(\frac{\log(\sqrt{D}/(pl))}{\log z}\Big).\end{align*}
Up to a minor contribution, the inner sum over $p$ can be relaxed to all primes $p\leqslant z.$ In fact, the terms with $p\mid \ell$ contribute at most
\begin{align*}
&\ll\frac{1}{(\log D)^4}\sum_{l\mid P(z)}\frac{1}{\beta^*(l)}m_2\Big(\frac{\log(\sqrt{D}/l)}{\log z}\Big)\sum_{p\mid l}\frac{\log p}{p}\\
&\ll\frac{1}{(\log D)^3\log \log D}\sum_{l\mid P(z)}\frac{1}{\beta^*(l)}m_2\Big(\frac{\log(\sqrt{D}/l)}{\log z}\Big)\\
&\ll\frac{1}{W(z)(\log D)^3\log\log D}.\end{align*}
We then derive that
\begin{align*}
S_{21}'(X)
&=-\frac{16(1+o(1))}{\cG(1)^2(\log D)^4}\sum_{l\mid P(z)}\frac{1}{\beta^*(l)}m_2\Big(\frac{\log(\sqrt{D}/l)}{\log z}\Big)\sum_{p<z}\frac{\log p}{\beta(p)}m_2\Big(\frac{\log(\sqrt{D}/(pl))}{\log z}\Big)\\
&\ \ \ \ \ \ +O\Big(\frac{\log z}{\log \log z}\frac{1}{W(z)(\log D)^4}\Big)\\
&=-\frac{32\tau c_{21}'(\tau)\log z}{\cG(1)^2W(z)(\log D)^4}\cdot(1+o(1)),\end{align*}
where
\begin{align*}
c_{21}'(\tau)=\int_0^1\int_0^1\sigma'((1-u)\tau)\ff(u\tau)\ff(u\tau-2v)\ud u\ud v.\end{align*}
In a similar manner, we can also show that
\begin{align*}
S_{21}''(X)&=\frac{32\tau c_{21}''(\tau)\log z}{\cG(1)^2W(z)(\log D)^4}\cdot(1+o(1)),\end{align*}
where
\begin{align*}
c_{21}''(\tau)=\int_0^1\int_0^1\sigma'((1-u)\tau)\ff(u\tau-2v)^2\ud u\ud v.\end{align*}
In conclusion, we obtain
\begin{align*}
S_{21}(X)=S_{21}'(X)+S_{21}''(X)
&=\frac{32\tau(c_{21}''(\tau)-c_{21}'(\tau))\log z}{\cG(1)^2W(z)(\log D)^4}\cdot(1+o(1)),\end{align*}

We now evaluate $S_{22}(X)$. For each squarefree $l\geqslant1$, we have
\begin{align*}
v(l)&=\beta(l)\sum_{u\mid l}\frac{\mu(u)}{\beta(u)}\sum_{p\mid l/u}\frac{\beta^*(p)\log p}{\beta(p)}\\
&=\beta(l)\sum_{p\mid l}\frac{\beta^*(p)\log p}{\beta(p)}\sum_{u\mid l/p}\frac{\mu(u)}{\beta(u)}\\
&=\beta(l)\sum_{p\mid l}\frac{\beta^*(l/p)\beta^*(p)\log p}{\beta(l/p)\beta(p)}\\
&=\beta^*(l)\log l.\end{align*}
Hence
\begin{align*}
S_{22}(X)&=\sum_{p<z}\beta^*(p)\log p\sum_{pl\mid P(z)}\beta^*(l)y_{pl}^2\\
&=\frac{16(1+o(1))}{\cG(1)^2(\log D)^4}\sum_{l\mid P(z)}\frac{1}{\beta^*(l)}\sum_{\substack{p<z\\p\nmid l}}\frac{\log p}{\beta^*(p)}m_2\Big(\frac{\log(\sqrt{D}/(pl))}{\log z}\Big)^2\end{align*}
by \eqref{eq:yl}. 
From partial summation, it follows that
\begin{align*}
S_{22}(X)&=\frac{32\tau c_{21}''(\tau)\log z}{\cG(1)^2W(z)(\log D)^4}\cdot(1+o(1)).\end{align*}

Combining all above evaluations, we find
\begin{align*}S(X,z;h,\boldsymbol\varrho)
&=\cG(1)X\{S_1(X)\cdot (\log X+c)-2S_{21}(X)+S_{22}(X)\}+O(X\cL^{-2})\\
&=(1+o(1))\frac{16\tau X\log z}{\cG(1)W(z)(\log D)^4}\Big\{c_1(\tau)\frac{\log X}{\log z}+4c_{21}'(\tau)-2c_{21}''(\tau))\Big\}.\end{align*}
Hence Theorem \ref{thm:asymptoticSelberg} follows by observing that $c_2(\tau)=2c_{21}'(\tau)-c_{21}''(\tau)$
and
\begin{align*}\cG(1)W(z)&=\prod_{p<z}\Big(1-\frac{1}{p}\Big)^2\cdot \prod_{p\geqslant z}\Big(1+\frac{h(p)}{p}\Big)\Big(1-\frac{1}{p}\Big)^2
=(1+o(1))\frac{\ue^{-2\gamma}}{(\log z)^2}\end{align*}
by Mertens' formula.

\bigskip

\section{Chebyshev approximation}
A lot of statistical analysis of $GL_2$ objects relies heavily on the properties of Chebychev polynomials $\{U_k(x)\}_{k\geqslant 0}$ with $x\in[-1,1],$ which can be defined recursively by
\[U_0(x)=1,\ \ U_1(x)=2x,\]
\[U_{k+1}(x)=2xU_k(x)-U_{k-1}(x),\ \ k\geqslant1.\]
It is well-known that Chebychev polynomials form an orthonormal basis of $L^2([-1, 1])$ with respect to the measure $\frac2{\pi}\sqrt{1-x^2}\ud x$. 
In fact, for any $f\in \mathcal{C}([-1,1])$, the expansion
\begin{align}\label{eq:psi-Chebyshevcoefficient}
f(x)=\sum_{k\geqslant 0}\beta_k(f)U_k(x)
\end{align}
holds with 
\[\beta_k(f)=\frac2{\pi} \int_{-1}^1f(t)U_k(t)\sqrt{1-t^2}\ud t.\] 
In practice, the following truncated approximation is usually more effective and useful, which has its prototype in \cite[Theorem 5.14]{MH03}.

\begin{lemma}\label{lm:Chebyshevapproximation}
Suppose $f:[-1,1]\rightarrow\bR$ has $C+1$ continuous derivatives on $[-1,1]$ with $C\geqslant2$. Then for each positive integer $K>C,$ there holds the approximation  
\[f(x)=\sum_{0\leqslant k\leqslant K}\beta_k(f)U_k(x)+O\Big(K^{1-C}\|f^{(C+1)}\|_1\Big)\]
uniformly in $x\in[-1,1]$, where the implied constant depends only on $C$.
\end{lemma}

\proof
For each $K>C$, we introduce the operator $\vartheta_{K}$ mapping $f \in \cC^{C+1}([-1,1])$ via  
\[(\vartheta_{K}f)(x):=\sum_{0\leqslant k\leqslant K}\beta_k(f)U_k(x)-f(x).\]
This gives the remainder of  approximation by Chebychev polynomials up to degree $K$.  Obviously, $(\vartheta_{K}f)(\cdot)\in \cC^{C+1}([-1,1])$ and in fact, $\vartheta_{K}$ is a bounded linear functional on $\cC^{C+1}([-1,1]),$ which vanishes on polynomials of degree $\leqslant K$. 

Using a theorem of Peano (Theorem 3.7.1, \cite{Da61}), we find that
\begin{align}\label{eq:theta_K}
(\vartheta_{K}f)(x)=\frac{1}{C!}\int_{-1}^1f^{(C+1)}(t)H_K(x,t)\ud t,
\end{align}
where
\begin{align*}H_K(x,t)=-\sum_{k>K}\lambda_k(t)U_k(x)\end{align*}
with
\begin{align*}\lambda_k(t)=\frac{2}{\pi}\int_t^1\sqrt{1-x^2}(x-t)^CU_k(x)\ud x.\end{align*}

Put $x=\cos\theta,t=\cos\phi$, so that
\begin{align*}\lambda_k(t) = \lambda_k(\cos \phi)=\frac{2}{\pi}\int_0^\phi(\cos\theta-\cos\phi)^C\sin\theta\sin((k+1)\theta)\ud\theta.\end{align*}
We deduce from integration by parts that
\begin{align*}\|\lambda_k\|_\infty\ll\frac{1}{k\dbinom{k-1}{C}},\end{align*}
where the implied constant is absolute.
For any $x,t\in[-1,1]$, the Stirling's formula $\log \Gamma(k) = (k-1/2)\log k - k +\log\sqrt{2\pi} + O(1/k)$ gives
\begin{align*}H_K(x,t)
&\ll\sum_{k>K}\frac{1}{\dbinom{k-1}{C}}=C!\sum_{k>K}\frac{\Gamma(k-C)}{\Gamma(k)}\\
& \ll \sum_{k>K} \Big(\frac{\ue}{k-C}\Big)^C \Big(1-\frac{C}k\Big)^{k-1/2}\\
& \ll K^{1-C}, 
\end{align*}
from which and \eqref{eq:theta_K} we conclude that
\begin{align*}
\|\vartheta_{K}f\|_\infty\ll K^{1-C}\|f^{(C+1)}\|_1.
\end{align*}
This completes the proof of the lemma.
\endproof

We now turn to derive a truncated approximation for $|x|$ on average.
\begin{lemma}\label{lm:|x|-Chebyshev}
Let $k,J$ be two positive integers and $K>1.$ Suppose $\{x_j\}_{1\leqslant j\leqslant J}\in[-1,1]$ and  $\by:=\{y_j\}_{1\leqslant j\leqslant J}\in\bC$ are two sequences satisfying
\begin{align}\label{eq:Chebyshev-averageassumption}
\max_{1\leqslant j\leqslant J}|y_j|\leqslant 1,\ \ \ \Bigg|\sum_{1\leqslant j\leqslant J}y_jU_k(x_j)\Bigg|\leqslant k^BU
\end{align}
with some $B\geqslant 1$ and $U>0$. Then we have
\begin{align*}
\sum_{1\leqslant j\leqslant J}y_j|x_j|
=\frac{4}{3\pi}\sum_{1\leqslant j\leqslant J}y_j+O\Big(UK^{B-1}(\log K)^{\delta(B)}+\frac{\|\mathbf{y}\|_1^2}{UK^B}\Big).
\end{align*}
where $\delta(B)$ vanishes unless $B=1$, in which case it is equal to $1$, and the $O$-constant depends only on $B.$
\end{lemma}

\proof
In order to apply Lemma \ref{lm:Chebyshevapproximation}, we would like to introduce a smooth function 
$R:[-1,1]\rightarrow[0,1]$ with $R(x)=R(-x)$
such that
\begin{align*}
\begin{cases}
R(x)=0,\ \ &x\in[-\varDelta,\varDelta],\\
R(x)=1,&x\in[-1,-2\varDelta]\cup[2\varDelta,1],\end{cases}
\end{align*}
where $\varDelta\in~]0,1[$ be a positive number to be fixed later.
We also assume the derivatives satisfy 
\begin{align*}
R^{(j)}(x)\ll_j \varDelta^{-j}
\end{align*}
for each $j\geqslant0$
with an implied constant depending only on $j$.

Put $f(x):=R(x)|x|.$ Due to smooth decay of $R$ at $x=0,$ we may apply Lemma \ref{lm:Chebyshevapproximation} to $f(x)$ with $C=2$,
getting
\begin{align*}
f(x)
&=\sum_{0\leqslant k\leqslant K}\beta_k(f)U_k(x)+O(K^{-1}\|f'''\|_1).
\end{align*}
Note that $f'''(x)$ vanishes unless  $x\in[-2\varDelta,-\varDelta]\cup[\varDelta,2\varDelta]$, in which case we have $f'''(x)\ll\varDelta^{-2}.$ It then follows that
\begin{align*}
f(x)
&=\sum_{0\leqslant k\leqslant K}\beta_k(f)U_k(x)+O\Big(\frac{1}{K\varDelta}\Big).
\end{align*}
Moreover, $f(x)-|x|$ vanishes unless $x\in[-2\varDelta,2\varDelta]$. This implies that $f(x)=|x|+O(\varDelta)$. In addition, $\beta_0(f)=\frac{4}{3\pi}+O(\varDelta)$. Therefore,
\begin{align*}
|x|&=\frac{4}{3\pi}+\sum_{1\leqslant k\leqslant K}\beta_k(f)U_k(x)+O\Big(\varDelta+\frac{1}{K\varDelta}\Big).
\end{align*}

We claim that 
\begin{align}\label{eq:beta_k(f)-upperbound}
\beta_k(f)\ll k^{-2}
\end{align}
for all $k\geqslant1$ with an absolute implied constant.
It then follows that
\begin{align*}
\sum_{1\leqslant j\leqslant J}y_j|x_j|-\frac{4}{3\pi}\sum_{1\leqslant j\leqslant J}y_j
&=\sum_{1\leqslant k\leqslant K}\beta_{k}(f)\sum_{1\leqslant j\leqslant J}y_jU_{k}(x_j)+O\Big(\|\mathbf{y}\|_1\varDelta+\frac{\|\mathbf{y}\|_1}{K\varDelta}\Big)\\
&\ll U\sum_{1\leqslant k\leqslant K}k^{B-2}+\|\mathbf{y}\|_1\varDelta+\frac{\|\mathbf{y}\|_1}{K\varDelta}\\
&\ll UK^{B-1}(\log K)^{\delta(B)}+\|\mathbf{y}\|_1\varDelta+\frac{\|\mathbf{y}\|_1}{K\varDelta},
\end{align*}
where the implied constant depends only on $B$.
To balance the first and last terms, we take $\varDelta=\|\mathbf{y}\|_1/(UK^B)$, which yields
\begin{align*}
\sum_{1\leqslant j\leqslant J}y_j|x_j|-\frac{4}{3\pi}\sum_{1\leqslant j\leqslant J}y_j
&\ll UK^{B-1}(\log K)^{\delta(B)}+\frac{\|\mathbf{y}\|_1^2}{UK^B}
\end{align*}
as expected.

It remains to prove the upper bound \eqref{eq:beta_k(f)-upperbound}. Since $U_k(\cos\theta)=\sin((k+1)\theta)/\sin\theta$, it suffices to show that
\begin{align}\label{eq:beta_k-upperbound}
\beta_k:=\int_0^{\frac{\pi}{2}}R(\cos\theta)(\sin2\theta)\sin((k+1)\theta)\ud\theta\ll k^{-2}
\end{align}
for all $k\geqslant3$ with an absolute implied constant.
From the elementary identity $2\sin\alpha\sin\beta=\cos(\alpha-\beta)-\cos(\alpha+\beta)$, it follows that
\begin{align*}
\beta_k=\int_0^{\arccos\varDelta}R(\cos\theta)(\sin2\theta)\sin((k+1)\theta)\ud\theta=\frac{\alpha(k-1,R)-\alpha(k+3,R)}{2},
\end{align*}
where, for $\ell\geqslant2$ and a function $g\in\cC^2([-1,1])$,
\begin{align*}
\alpha(\ell,g):=\int_0^{\arccos\varDelta}g(\cos\theta)\cos(\ell\theta)\ud\theta.
\end{align*}

From integration by parts, we derive that
\begin{align*}
\alpha(\ell,g)&=\frac{1}{\ell}\int_0^{\arccos\varDelta}g'(\cos\theta)(\sin\theta)\sin(\ell\theta)\ud\theta=\frac{\alpha(\ell-1,g')-\alpha(\ell+1,g')}{2\ell},
\end{align*}
and also
\begin{align*}
\alpha(\ell,g')&=\frac{\alpha(\ell-1,g'')-\alpha(\ell+1,g'')}{2\ell}.
\end{align*}
It then follows that
\begin{align*}
\alpha(\ell,g)&=\frac{\alpha(\ell-2,g'')-\alpha(\ell,g'')}{4\ell(\ell-1)}-\frac{\alpha(\ell,g'')-\alpha(\ell+2,g'')}{4\ell(\ell+1)}.
\end{align*}
We then further have
\begin{align*}
\beta_k&=\frac{1}{8}(\beta_{k,1}-\beta_{k,2})
\end{align*}
with
\begin{align*}
\beta_{k,1}&=\frac{\alpha(k-3,R'')-\alpha(k-1,R'')}{(k-1)(k-2)}-\frac{\alpha(k-1,R'')-\alpha(k+1,R'')}{k(k-1)},\\
\beta_{k,2}&=\frac{\alpha(k-1,R'')-\alpha(k+1,R'')}{k(k+1)}-\frac{\alpha(k+1,R'')-\alpha(k+3,R'')}{(k+1)(k+2)}.
\end{align*}
Note that
\begin{align*}
\alpha(k-3,R'')-\alpha(k-1,R'')&=\int_{\arccos2\varDelta}^{\arccos\varDelta}R''(\cos\theta)\{\cos((k-3)\theta)-\cos((k-1)\theta)\}\ud\theta\\
&=2\int_{\arccos2\varDelta}^{\arccos\varDelta}R''(\cos\theta)(\sin(k-2)\theta)(\sin\theta)\ud\theta
\end{align*}
and
\begin{align*}
\alpha(k-1,R'')-\alpha(k+1,R'')
&=2\int_{\arccos2\varDelta}^{\arccos\varDelta}R''(\cos\theta)(\sin k\theta)(\sin\theta)\ud\theta.
\end{align*}
Hence
\begin{align*}
\beta_{k,1}
&=\frac{2}{(k-1)(k-2)}\int_{\arccos2\varDelta}^{\arccos\varDelta}R''(\cos\theta)(\sin(k-2)\theta)(\sin\theta)\ud\theta\\
&\ \ \ \ \ -\frac{2}{k(k-1)}\int_{\arccos2\varDelta}^{\arccos\varDelta}R''(\cos\theta)(\sin k\theta)(\sin\theta)\ud\theta\\
&=\frac{2}{(k-1)(k-2)}\int_{\arccos2\varDelta}^{\arccos\varDelta}R''(\cos\theta)\{\sin(k-2)\theta-\sin k\theta\}(\sin\theta)\ud\theta\\
&\ \ \ \ \ +\frac{4}{k(k-1)(k-2)}\int_{\arccos2\varDelta}^{\arccos\varDelta}R''(\cos\theta)(\sin k\theta)(\sin\theta)\ud\theta.
\end{align*}
The first term can be evaluated as
\begin{align*}
&=\frac{2}{(k-1)(k-2)}\int_{\arccos2\varDelta}^{\arccos\varDelta}R''(\cos\theta)(\sin(k-1)\theta)(\cos\theta)(\sin\theta)\ud\theta\\
&\ll\frac{1}{k^2}\int_{\arccos2\varDelta}^{\arccos\varDelta}\varDelta^{-2}\cos\theta\ud\theta\ll\frac{1}{k^2}.
\end{align*}
Again from the integration by parts, the second term is
\begin{align*}
=\frac{4}{(k-1)(k-2)}\int_{\arccos2\varDelta}^{\arccos\varDelta}R'(\cos\theta)(\cos k\theta)\ud\theta
\ll\frac{1}{k^2}\int_{\arccos2\varDelta}^{\arccos\varDelta}\varDelta^{-1}\ud\theta\ll\frac{1}{k^2}.
\end{align*}
Hence $\beta_{k,1}\ll k^{-2}$, and similarly $\beta_{k,2}\ll k^{-2}.$
These yield \eqref{eq:beta_k-upperbound}, and thus \eqref{eq:beta_k(f)-upperbound}, which completes the proof of the lemma.
\endproof

Note that $U_k(\cos\theta)=\sym_k(\theta).$ Taking $x_j=\cos\theta_j$ in Lemma \ref{lm:|x|-Chebyshev}, we obtain the following truncated approximation for $|\cos|$.
\begin{lemma}\label{lm:cos-Chebyshev}
Let $k,J$ be two positive integers and $K>1.$ Suppose $\{\theta_j\}_{1\leqslant j\leqslant J}\in[0,\pi]$ and  $\by:=\{y_j\}_{1\leqslant j\leqslant J}\in\bC$ are two sequences satisfying
\begin{align*}
\max_{1\leqslant j\leqslant J}|y_j|\leqslant 1,\ \ \ \Bigg|\sum_{1\leqslant j\leqslant J}y_j\sym_k(\theta_j)\Bigg|\leqslant k^BU
\end{align*}
with some $B\geqslant 1$ and $U>0.$ Then we have
\begin{align*}
\sum_{1\leqslant j\leqslant J}y_j|\cos\theta_j|
=\frac{4}{3\pi}\sum_{1\leqslant j\leqslant J}y_j+O\Big(UK^{B-1}(\log K)^{\delta(B)}+\frac{\|\mathbf{y}\|_1^2}{UK^B}\Big).
\end{align*}
where $\delta(B)$ is defined as in Lemma {\rm \ref{lm:|x|-Chebyshev}.}
\end{lemma}

\smallskip

\bibliographystyle{plainnat}

\bigskip

\end{document}